\documentclass[12pt,a4paper,reqno]{amsart}
\usepackage{amssymb}
\usepackage{dcpic,pictex}
\usepackage{multirow}
\usepackage{graphicx}
\usepackage{slashbox}
\usepackage{cite}

\usepackage{hyperref}
\newtheorem{Def}{Definition}[section]
\newtheorem{Ex}{Example}[section]
\newtheorem{The}{Theorem}[section]
\newtheorem{Rem}{Remark}[section]

\begin{document}


\title[Subordination principles for the fractional diffusion-wave equation]{Subordination principles for the multi-dimensional space-time-fractional diffusion-wave equation}

\author{Yuri Luchko}
\curraddr{Beuth Technical University of Applied Sciences Berlin,  
     Department of  Mathematics, Physics, Chemistry,  
     Luxemburger Str. 10,  
     13353 Berlin,
Germany}
\email{luchko@beuth-hochschule.de}

\date{13.02.2018}
\dedicatory{}
\thanks{}

\begin{abstract}
This paper is devoted to an in deep investigation of the first fundamental solution to the linear multi-dimensional space-time-fractional diffusion-wave equation. This equation is obtained from the diffusion equation by replacing the first order time-deri\-va\-ti\-ve by  the Caputo fractional derivative of order $\beta,\ 0 <\beta \leq 2$ and the Laplace operator by the fractional Laplacian $(-\Delta)^{\frac\alpha 2}$ with $0<\alpha \leq 2$. First, a representation of the fundamental solution in form of a Mellin-Barnes integral is deduced  by employing the technique of the Mellin integral transform. This representation is then used for establishing of several subordination formulas that connect the fundamental solutions for different values of the fractional derivatives $\alpha$ and $\beta$.  We also discuss some new cases of completely monotone functions and probability density functions that are expressed in terms of the Mittag-Leffler function, the Wright function, and the generalized Wright function. 
\end{abstract}



\maketitle

\section{Introduction}

A subordination principle for completely positive measures was discussed in detail in \cite{Pru93} and applied there for constructing new resolvents for the abstract Volterra integral equations based on the known ones. In \cite{Baz00} and \cite{Baz01}, this subordination principle was extended and specialized for the abstract fractional evolution equations in the form
\begin{equation}
\label{1.1}
D^\beta u(t) = Au(t),
\end{equation}
subject to the initial conditions
\begin{equation}
\label{1.2}
u(0) = x,\ u^{(k)}(0) = 0,\ \, k=0,\dots, n-1,
\end{equation}
where $D^\beta$ is the Caputo fractional derivative of order $\beta$ that will be defined in the next section, $n-1 <\beta \le n,\ n \in \mathbb{N}$, and $A$ is a linear closed unbounded operator densely defined in a Banach space $X$, where the initial condition from \eqref{1.2} belongs to, i.e., $x \in X$. 

Let $S_\beta(t)$ be a solution operator to the abstract initial-value problem \eqref{1.1}-\eqref{1.2}, $0 <\beta <\delta \le 2$, $\gamma = \beta/\delta$. Then the subordination formula 
\begin{equation}
\label{1.3}
S_\beta(t)x = \int_0^\infty t^{-\gamma} W_{1-\gamma,-\gamma}(-st^{-\gamma})S_\delta(s)x \ ds,\ t>0,\ x\in X
\end{equation}
is valid under some conditions on the operator $A$ (see \cite{Baz00} and \cite{Baz01} for details). The function $W_{1-\gamma,-\gamma}(-\tau)$ from \eqref{1.3} is a special case of the Wright function that will be introduced in the next section. It is important to mention that this function is non-negative for $\tau \in \mathbb{R_+}$ and can be interpreted as a probability density function.

Very recently, the subordination principle was extended to the case of the multi-term time-fractional diffusion-wave equations in \cite{Baz17} and to the case of the distributed order time-factional evolution equations in the Caputo and Riemann-Liouville sense in \cite{Baz15}. 

All publications mentioned above deal with the abstract evolution equations with a  linear closed unbounded operator $A$ subject to some additional conditions. In Fractional Calculus and its applications, an important particular case of these equations, namely,  fractional differential equations with both time-factional and space-fractional derivatives are nowadays subject of very intensive research. Say, in \cite{Han2}, mathematical, physical, and probabilistic properties of the fundamental solutions to the multi-dimensional  space-time-fractional diffusion-wave equation were considered. In the papers \cite{BL1}, \cite{BL2}, \cite{L2013}-\cite{L2017} the method of the Mellin-Barnes integral representations was employed to derive further properties of solutions to the multi-dimensional space-time-fractional diffusion-wave equation and its important particular cases  as the $\alpha$-fractional diffusion equations and the $\alpha$-fractional wave equation. Still, the scope of the properties, particular cases, integral and series representations, asymptotic formulas, etc. known for the fundamental solution to the one-dimensional space-time-fractional diffusion-wave equation (see \cite{MLP2001} for its detailed theory) is essentially more expanded compared to the multi-dimensional case and thus further investigations of the multi-dimensional case are required. 

In this paper, both known and new subordination formulas for the fundamental solutions  to the Cauchy problems for the multi-dimensional space-time-fractional diffusion-wave equation are derived and discussed. The subordination formulas that connect the fundamental solutions for different orders of the time-fractional derivative of the type given by \eqref{1.3} are already known, but here we apply a different method for their derivation. To the best of the authors knowledge, the  subordination formulas presented in this paper that connect the fundamental solution to the Cauchy problem for the multi-dimensional space-time-fractional diffusion-wave equation with the fundamental solution of the conventional diffusion equation as well as a subordination formula for the space-fractional diffusion equation are new. For the subordination formulas for the fundamental solutions to the one-dimensional space-time fractional diffusion-wave equation we refer to \cite{MLP2001}.

The rest of the paper is organized as follows. In the second section, we formulate the problem we deal with and recall the Mellin-Barnes representations of the fundamental solution to the Cauchy problem for the multi-dimensional space-time-fractional diffusion-wave equation that were derived in the previous publications of the author and his co-authors. The third section is devoted to a  discussion of a special technique for derivation of new completely monotone functions and new non-negative functions that can be interpreted as probability density functions. This technique is then applied  for construction of some new completely monotone functions and probability density functions in terms of the Mittag-Leffler function, the Wright function, and the generalized Wright function that will be used in the further discussions. In the final section of the paper, the Mellin-Barnes representations of the fundamental solution are employed for derivation of both known and new subordination formulas for the solutions  to the Cauchy problem for the multi-dimensional space-time-fractional diffusion-wave equation with different orders of the time- and space-fractional derivatives. 

\section{Problem formulation and auxiliary results}

In this paper, we deal with the linear multi-dimensional space-time-fractional diffusion-wave equation in the following form:
\begin{equation}\label{eq1}
    D_t^{\beta}u(\mathrm{x},t)=-(-\Delta)^{\frac\alpha 2}u(\mathrm{x},t), \quad \mathrm{x}\in
\mathbb{R}^n,\; t>0,\; 0<\alpha \leq 2,\; 0 < \beta \leq 2.
\end{equation}
In equation \eqref{eq1}, $D_t^{\beta}$ denotes the Caputo time-fractional derivative of  order $\beta,\ \beta >0$ defined by the formula 
\begin{equation}\label{eq2}
    D_t^\beta u(\mathrm{x},t)=\left(I^{n-\beta}_t \frac{\partial^nu}{\partial t^n}\right)(t), \quad n-1<\beta\leq n,\ n\in\mathbb{N}\, ,
\end{equation}
where $I^{\gamma}_t$ is the Riemann-Liouville fractional integral: 
$$
(I^\gamma_t u)(t)= 
\begin{cases}
\frac 1{\Gamma(\gamma)}\int_0^t (t-\tau)^{\gamma-1}u(\mathrm{x},\tau)\, d\tau\ \mbox{ for }\ \gamma> 0, \\
u(\mathrm{x},t)\ \mbox{ for }\ \gamma= 0. 
\end{cases}
$$
The fractional Laplacian $(-\Delta)^{\frac\alpha 2}$ from the equation \eqref{eq1} is defined as a pseudo-differential operator with the symbol $|\kappa|^\alpha$ (\cite{SZ1997,SKM1993}):  
\begin{equation}\label{eq3}
    \left( {\mathcal F}\, (-\Delta)^{\frac{\alpha}{2}}f\right) (\kappa)=|\kappa|^\alpha({\mathcal F}\, f)(\kappa)\, ,
\end{equation}
where $({\mathcal F}\, f)(\kappa)$ is the Fourier transform of a function $f$ at the point $\kappa \in\mathbb{R}^n$ defined by
\begin{equation}\label{eq4}
    ({\mathcal F}\, f)(\kappa)=\hat{f}(\kappa)=\int_{\mathbb{R}^n} e^{i\kappa \cdot \mathrm{x}}f(\mathrm{x})\, d\mathrm{x}\, .
\end{equation}
For $0<\alpha<m,\ m\in \mathbb{N}$ and $x\in \mathbb{R}^n$,  the fractional Laplacian can be also represented as a hypersingular integral (\cite{SKM1993}):
\begin{equation}
\label{hyp}
(-\Delta)^{\frac{\alpha}{2}}f(\mathrm{x})=\frac 1{d_{n,m}(\alpha)} \int_{\mathbb{R}^n}\frac{\left(\Delta^m_{\mathrm{h}} f\right)(\mathrm{x})}{|\mathrm{h}|^{n+\alpha}}\, d\mathrm{h} 
\end{equation}
with a suitably defined finite differences operator $\left(\Delta^m_{\mathrm{h}} f\right)(\mathrm{x})$ and 
a normalization constant $d_{n,m}(\alpha)$. 

The representation \eqref{hyp} of the fractional Laplacian in form of the hypersingular integral does not depend on $m, \ m\in \mathbb{N}$ provided $\alpha<m$ (\cite{SKM1993}). 

In the one-dimensional case, the equation \eqref{eq1} is a particular case of a more general equation with the Caputo time-fractional derivative and the Riesz-Feller space-fractional derivative that was discussed in detail in \cite{MLP2001}. For $\alpha = 2$, the fractional Laplacian $(-\Delta)^{\frac\alpha 2}$ is just $-\Delta$ and thus the equation \eqref{eq1} is a particular case of the time-fractional diffusion-wave equation that was considered in many publications including, say, \cite{EK2004}, \cite{FV2016}, \cite{Han3}, \cite{K1990}, \cite{K2014}, \cite{L2010}, and \cite{SW1989}. For $\alpha =2$ and $\beta = 1$, the equation  \eqref{eq1} is reduced to the diffusion equation and for  $\alpha =2$ and $\beta = 2$ it is the wave equation that justifies its denotation as a fractional diffusion-wave equation. 

In this paper, we consider the Cauchy problem for the space-time-fractional diffusion-wave equation  \eqref{eq1} with the Dirichlet initial conditions: 
\begin{equation}
\label{eq6}
    u(\mathrm{x},0)=\varphi(\mathrm{x})\, , \quad \mathrm{x}\in \mathbb{R}^n  
\end{equation}
if the order $\beta$ of the time-fractional derivative satisfies the condition $0 < \beta \leq 1$ or  
\begin{equation}\label{eq6_1}
    u(\mathrm{x},0)=\varphi(\mathrm{x})\, ,\ \frac{\partial u}{\partial t}(\mathrm{x},0)=0\, , \quad \mathrm{x}\in \mathbb{R}^n  
\end{equation}
if $1 < \beta \leq 2$.

Because the initial-value problem \eqref{eq1}, \eqref{eq6} or  \eqref{eq1}, \eqref{eq6_1}, respectively, is a linear one, its solution can be 
represented in the form
\begin{equation}
\label{FS}
u(\mathrm{x},t) = \int_{\mathbb{R}^n} G_{\alpha,\beta, n}(\zeta,t)\varphi(\mathrm{x}-\zeta)\, d\zeta,
\end{equation}
where $G_{\alpha,\beta,n}$ is  the so-called first fundamental solution to the fractional diffusion-wave equation \eqref{eq1} and the function $\varphi$ is given in the initial condition. By $G_{\alpha,\beta,n}$, the solution to the equation \eqref{eq1} with the initial condition  ($0 < \beta \leq 1$)
$$
u(\mathrm{x},0)=\prod_{i=1}^n\delta(x_i)\, , \quad \mathrm{x}=(x_1,x_2,\ldots,x_n)\in \mathbb{R}^n
$$
or the initial conditions ($1 < \beta \leq 2$)
$$
u(\mathrm{x},0)=\prod_{i=1}^n\delta(x_i),\ \frac{\partial u}{\partial t}(\mathrm{x},0)=0\, , \quad \mathrm{x}=(x_1,x_2,\ldots,x_n)\in \mathbb{R}^n,
$$
respectively, is denoted
with $\delta$ being the Dirac delta function.

Thus the behavior of the solutions to the problem \eqref{eq1}, \eqref{eq6} or  \eqref{eq1}, \eqref{eq6_1}, respectively, is determined by the fundamental solution $G_{\alpha,\beta,n}$ and the focus of this paper is on derivation  of some new properties of the fundamental solution. In particular, we deal with the subordination formulas for the fundamental solution $G_{\alpha,\beta,n}$ in the form 
\begin{equation}
\label{FS-S1}
G_{\alpha,\beta,n}(\mathrm{x},t) = \int_0^\infty \Phi(s,t) G_{\hat{\alpha},\hat{\beta},n}(\mathrm{x},s)\, ds,
\end{equation}
where the kernel function $\Phi=  \Phi(s,t)$ can be interpreted as a probability density function in $s,\ s\in \mathbb{R}_+$ for each value of $t,\ t>0$. Let us note here that any subordination formula for the solution operator $S_{\alpha,\beta,n}(t)$ to the initial-value problem \eqref{eq1}, \eqref{eq6} or  \eqref{eq1}, \eqref{eq6_1}, respectively, in the form (see e.g. \eqref{1.1})
\begin{equation}
\label{FS-S2}
S_{\alpha,\beta,n}(t)\varphi = \int_0^\infty \Phi(s,t) S_{\hat{\alpha},\hat{\beta},n}(s)\varphi\, ds
\end{equation}
induces a subordination formula of the type \eqref{FS-S1} for the fundamental solution $G_{\alpha,\beta,n}$ just by setting $\varphi$ to be the Dirac $\delta$-function. Vice versa, any subordination formula of type \eqref{FS-S1} for the fundamental solution $G_{\alpha,\beta,n}$ automatically leads to a subordination formula for the solution operator $S_{\alpha,\beta,n}(t)$ of type \eqref{FS-S2} because of the representation \eqref{FS}. Indeed, we have the following chain of (formal) transformations:
$$
S_{\alpha,\beta,n}(t)\varphi = \int_{\mathbb{R}^n} G_{\alpha,\beta, n}(\zeta,t)\varphi(\mathrm{x}-\zeta)\, d\zeta \ = \ 
\int_{\mathbb{R}^n} \int_0^\infty \Phi(s,t) G_{\hat{\alpha},\hat{\beta},n}(\zeta,s)\, ds\, \varphi(\mathrm{x}-\zeta)\, d\zeta =
$$
$$
\int_0^\infty   \Phi(s,t) \int_{\mathbb{R}^n} G_{\hat{\alpha},\hat{\beta},n}(\zeta,s)\, \varphi(\mathrm{x}-\zeta)\, d\zeta \, ds \ = \  \int_0^\infty \Phi(s,t) S_{\hat{\alpha},\hat{\beta},n}(s)\varphi\, ds.
$$
Thus a (formal) derivation of the subordination formulas for the solution operator  can be reduced to derivation of the subordination formulas for the fundamental solution. Of course, afterwards, the subordination formulas for the solution operator $S_{\alpha,\beta,n}(t)$ should be strictly proved. In this paper, we restrict ourselves to the first step of this procedure, namely, to derivation of some subordination formulas for the fundamental solution $G_{\alpha,\beta,n}$. Their translation to the solution operator $S_{\alpha,\beta,n}(t)$ will be considered elsewhere. 

The subordination formulas for the fundamental solution will be deduced based on  their   Mellin-Barnes representations. For the reader's convenience, a short sketch of derivation of these representations will be presented in the rest of this section. For the details we refer to \cite{L2014}, \cite{L2015} for the case $\beta= \alpha$, to \cite{BL2} for the case $\beta= \alpha/2$, and to \cite{BL1}, \cite{L2017}  for the general case. 

Application of the multi-dimensional Fourier transform \eqref{eq4} with respect to the spatial variable $\mathrm{x}\in \mathbb{R}^n$ to the equation \eqref{eq1} and to the initial conditions \eqref{eq6} or \eqref{eq6_1}, respectively, with $\varphi(\mathrm{x})=\prod_{i=1}^n\delta(x_i)$ leads to the ordinary fractional differential equation  in the Fourier domain
\begin{equation}
\label{trans}
D_t^{\beta}\hat{G}_{\alpha,\beta,n}(\kappa,t)+|\kappa|^\alpha\hat{G}_{\alpha,\beta,n}(\kappa,t)=0,  
\end{equation}
along with the initial conditions
\begin{equation}
\label{trans_1}
\hat{G}_{\alpha,\beta,n}(\kappa,0)=1
\end{equation}
in the case $0 <\beta \leq 1$ or with the initial conditions
\begin{equation}
\label{trans_2}
\hat{G}_{\alpha,\beta,n}(\kappa,0)=1,\ \frac{\partial }{\partial t}\hat{G}_{\alpha,\beta,n}(\kappa,0)=0
\end{equation}
in the case $1 <\beta \leq 2$.  

In both cases, the unique solution of \eqref{trans} with the initial conditions \eqref{trans_1} or  \eqref{trans_2}, respectively,    has the following form (see e.g. \cite{L1999}):
\begin{equation}
\label{eq7}
    \hat{G}_{\alpha,\beta,n}(\kappa,t)=E_{\beta}\left(-|\kappa|^\alpha t^{\beta}\right)\, 
\end{equation}
in terms of  the Mittag-Leffler function $E_{\beta}(z)$ that is defined by a convergent series 
\begin{equation}
\label{eq8}
    E_{\beta}(z)=\sum_{n=0}^\infty \frac{z^n}{\Gamma(1+\beta\, n)}\, , \quad \beta>0,\ z\in \mathbb{C}.
\end{equation}
Under the condition 
$\alpha>1$, 
one has the inclusion $\hat{G}_{\alpha,\beta,n} \in L_1(\mathbb{R}^n)$ because of the asymptotic formula (see e.g. \cite{E1955})
\begin{equation}
\label{M-L-a}
E_\beta(-x) = -\sum_{k=1}^{m} \frac{(-x)^{-k}}{\Gamma(1-\beta k)} \ 
+O(|x|^{-1-m}),\ m\in \mathbb{N},\ x\to +\infty,\ 0<\beta <2.
\end{equation}
Thus the inverse Fourier transform of 
\eqref{eq7} can be represented as follows
\begin{equation}\label{eq9}
    G_{\alpha,\beta,n}(\mathrm{x},t)=\frac 1{(2\pi)^n}\int_{\mathbb{R}^n} e^{-i\kappa \cdot \mathrm{x}}E_{\beta}\left(-|\kappa|^\alpha t^{\beta}\right)\, d\kappa\, , \quad\mathrm{x}\in\mathbb{R}^n\, , t>0\, .
\end{equation}
Because  $E_{\beta}\left(-|\kappa|^\alpha t^{\beta}\right)$ is a radial function, the known formula (see e.g. \cite{SKM1993})
\begin{equation}\label{eq10}
    \frac 1{(2\pi)^n}\int_{\mathbb{R}^n} e^{-i\kappa \cdot \mathrm{x}}\varphi(|\kappa|) \, d\kappa= \frac{|\mathrm{x}|^{1-\frac{n}{2}}}{(2\pi)^{\frac n2}}\, \int_0^\infty \varphi (\tau) \tau^{\frac n2} J_{\frac n2-1}(\tau|\mathrm{x}|) \, d\tau\, 
\end{equation}
for the Fourier transform of the radial functions can be applied, where $J_\nu$ denotes the Bessel function with the index $\nu$ (for the properties of the the Bessel function see e.g. \cite{E1953}), and we arrive at the representation
\begin{equation}\label{eq11}
    G_{\alpha,\beta, n}(\mathrm{x},t)=\frac{|\mathrm{x}|^{1-\frac n 2}}{(2\pi)^{\frac n 2}}\, \int_0^\infty E_{\beta}\left(-\tau^\alpha t^{\beta}\right)  \tau^{\frac n 2} J_{\frac n 2-1}(\tau|\mathrm{x}|) \, d\tau\, ,
\end{equation}
whenever the integral in \eqref{eq11} converges absolutely or at least conditionally.

The representation \eqref{eq11} can be transformed to a Mellin-Barnes integral. 
We start with the case $|\mathrm{x}|=0$ ($\mathrm{x}=(0, \ldots, 0)$) and get the formula
$$
G_{\alpha,\beta,n}(0,t)=\frac{1}{(2\pi)^{n}}\int_{\mathbb{R}^{n}} E_{\beta} (-|\kappa|^{\alpha}t^{\beta})d\kappa
$$
that can be represented in the form
\begin{equation}\label{eq13}
G_{\alpha,\beta,n}(\mathrm{0},t)=\frac{1}{(2\pi)^{n}}\frac{2\pi^{\frac{n}{2}}}{\Gamma(\frac{n}{2})}\int_0^{\infty} E_{\beta}(-\tau^{\alpha}t^\beta)\, \tau^{n-1}\, d\tau
\end{equation}
due to the known formula (see  e.g. \cite{SKM1993})
\begin{equation}\label{eq12}
\int_{\mathbb{R}^{n}} f(|\mathrm{x}|)d\mathrm{x}=\frac{2\pi^{\frac{n}{2}}}{\Gamma(\frac{n}{2})}\int_0^{\infty} \tau^{n-1}f(\tau)d\tau.
\end{equation}
The asymptotics of 
the Mittag-Leffler function 
ensures convergence of the integral in \eqref{eq13} under the 
condition $0<n<\alpha$. Say, for  
$1< \alpha \leq 2$ the fundamental solution  $G_{\alpha, \beta,n}$ is finite at $|x|=0$ only in the one-dimensional case. In this case, we get the formula 
$$
G_{\alpha,\beta,1}(\mathrm{0},t)=\frac{t^{-\frac{\beta}{\alpha}}}{\alpha \pi}\int_0^{\infty} E_{\beta}(-u)\, u^{\frac{1}{\alpha}-1}\, du\ = \frac{t^{-\frac{\beta}{\alpha}}}{\alpha\pi} \frac{\Gamma\left(\frac{1}{\alpha}\right)\Gamma\left(1-\frac{1}{\alpha}\right)}
{\Gamma\left( 1- \frac{\beta}{\alpha}\right)}
$$
that is valid for $\alpha >1$ if $0 < \beta < 2$ and for $\alpha >2$ if $\beta = 2$. 
This formula is nothing else as an easy consequence from the known Mellin integral transform of the Mittag-Leffler function (see e.g. \cite{LK2013}, \cite{Mar83}):
\begin{equation}
\label{ML-Mellin1}
\int_0^{\infty} E_{\beta}(-t)\, t^{s-1}\, dt \ = \ 
\frac{\Gamma(s)\Gamma(1-s)}
{\Gamma(1 -\beta s)}\ \mbox{ if } \ 
\begin{cases}
0 < \Re(s) < 1\ \mbox{ for } \  0<\beta < 2, \\  
0 < \Re(s) < 1/2\ \mbox{ for } \  \beta=2.
\end{cases}
\end{equation}
If the dimension $n$ of the equation \eqref{eq1} is greater that one, the fundamental solution $G_{\alpha, \beta,n}(x,t)$ has an integrable singularity at the point $|x|=0$.

The Mellin integral transform plays an important role in Fractional Calculus in general and for derivation of the results of this paper in particular, so let us recall the definitions of the Mellin transform and the inverse Mellin transform, respectively: 
\begin{equation}\label{eq15}
f^\ast(s)=({\mathcal M} f(t))(s)=\int_0^{\infty} f(t) t^{s-1}  \, dt\, , \quad t>0\, ,
\end{equation}
\begin{equation}\label{eq16}
f(t)=({\mathcal M}^{-1}f^\ast(s))(t)=\frac 1{2\pi i}\int_{\gamma-i\infty}^{\gamma+i\infty} f^\ast(s)t^{-s} \, ds\, ,  \gamma_1<\Re(s)=\gamma<\gamma_2\, .
\end{equation}
As it is well known, the Mellin integral transform exists for the functions continuous on the intervals $(0,\varepsilon]$ and $[E,+\infty)$ and integrable on  the interval $(\epsilon,E)$ with any $\varepsilon,\ E,$ $0<\varepsilon<E<+\infty$  that satisfy the estimates $|f(t)|\leq M_1t^{-\gamma_1}$ for $0<t<\varepsilon$ and $|f(t)|\leq M_2t^{-\gamma_2}$ for $t>E$ with $\gamma_1<\gamma_2$ and some constants $M_1,\ M_2$. In this case, the Mellin integral transform $f^\ast(s)$  is analytic in the vertical strip $\gamma_1<\Re(s)=\gamma<\gamma_2$.

If $f$ is piecewise differentiable  and $t^{\gamma-1}f(t)\in L^c(0,\infty)$, then the inversion formula \eqref{eq16} holds at all points of continuity of the function $f$. The integral in the formula \eqref{eq16} has to be considered in the sense of the Cauchy principal value.

For the general theory of the Mellin integral transform we refer the reader to \cite{Mar83}. Several applications of the Mellin integral transform in fractional calculus are discussed in  \cite{LK2013}. 

In the further discussions, we employ some of the elementary
properties of the Mellin
integral transform that are summarized below. Denoting by   $ \overset{\mathcal M}{\longleftrightarrow} $
the juxtaposition of a function $f$ with its
Mellin transform $f^*\,,$ the needed rules are:
 \begin{eqnarray}
\label{(1.51)}
f(at) & \overset{\mathcal M}{\longleftrightarrow}  & a^{-s}f^*(s), \ a>0, \\
\label{(1.52)}
t^\alpha f(t) & \overset{\mathcal M}{\longleftrightarrow}  & f^*(s+\alpha), \\
\label{(1.53)}
  f(t^\alpha) & \overset{\mathcal M}{\longleftrightarrow}  & \frac{1}{\vert \alpha \vert} f^*(s/\alpha),\  \alpha\not=0.
\end{eqnarray}
Another important operational relation is the convolution theorem for the Mellin integral transform that reads as follows:
\begin{equation}
\label{eq17}
\int_{0}^{\infty} f_1(\tau)f_2\left(\frac y \tau\right)\, \frac{d \tau}{\tau} \ \overset{\mathcal M}{\longleftrightarrow} \  f_1^\ast(s)f_2^\ast(s).
\end{equation}

Now we proceed with the case $\mathrm{x}\not = 0$. As to the convergence of the integral in \eqref{eq11}, it follows from the asymptotic formulas for the Mittag-Leffler function and the known asymptotic behavior of the Bessel function (see e.g. \cite{E1953}) that it converges conditionally in the case $n<2\alpha +1$ and absolute in the case  $n<2\alpha-1$. Say, for $1<\alpha \leq 2$  and $n=1,2,3$ the integral in \eqref{eq11} is at least conditionally convergent.

It can be easily seen that for $\mathrm{x}\neq 0$ the integral at the right-hand side of the formula \eqref{eq11} is nothing else as the Mellin convolution of the functions
$$
f_1(\tau)=E_{\beta}(-\tau^{\alpha}\, t^{\beta} ) \quad \mbox{and}\quad 
f_2(\tau)=\frac{|\mathrm{x}|^{-n}}{(2\pi)^{\frac n 2}}\,  
\tau^{-\frac{n}{2}-1}\, J_{\frac{n}{2}-1}
\left(\frac{1}{\tau}\right)
$$
at the point $y=\frac{1}{|\mathrm{x}|}$.

The Mellin transform of the Mittag-Leffler function \eqref{ML-Mellin1}, the known Mellin integral transform of the Bessel function (\cite{Mar83})
$$
J_\nu (2\sqrt{\tau}) \overset{\mathcal M}{\longleftrightarrow} \frac{\Gamma(\nu/2 +s)}{\Gamma(\nu/2 +1 -s)},\ -\Re(\nu/2) < \Re(s) < 3/4,
$$
and some elementary properties of the Mellin integral transform (see e.g. \cite{LK2013,Mar83}) lead to the Mellin transform formulas:
\begin{eqnarray*}
f_1^{\ast}(s)=\frac{t^{-\frac{\beta}{\alpha}\, s}}{\alpha} \frac{\Gamma(\frac{s}{\alpha})\Gamma(1-\frac{s}{\alpha})}{\Gamma(1-\frac{\beta}{\alpha}\, s)}\, , \quad 0<\Re(s)<\alpha\, ,\\
f_2^{\ast}(s)=\frac{|\mathrm{x}|^{-n}}{(2\pi)^{\frac n 2}} 
\left(\frac{1}{2}\right)^{-\frac{n}{2}+s} \frac{\Gamma\left(\frac{n}{2}-\frac{s}{2}\right)}
{\Gamma\left(\frac{s}{2}\right)}\, , \quad \frac{n}{2} - \frac{1}{2}<\Re(s)<n\, .
\end{eqnarray*}
These two formulas, the convolution theorem  \eqref{eq17} for the Mellin transform,  and the inverse Mellin transform formula \eqref{eq16} result in the following Mellin-Barnes integral   representation
of the fundamental solution $G_{\alpha,\beta,n}$:
\begin{equation}
\label{eq20}
G_{\alpha,\beta,n}(\mathrm{x},t)=\frac{1}{\alpha}\frac{|\mathrm{x}|^{-n}}{\pi^{\frac{n}{2}}}\frac{1}{2\pi i} \int_{\gamma-i\infty}^{\gamma+i\infty} \frac{\Gamma\left(\frac{n}{2}-\frac{s}{2}\right)\Gamma\left(\frac{s}{\alpha}\right) \Gamma\left(1-\frac{s}{\alpha}\right)}{\Gamma\left(1-\frac{\beta}{\alpha}s\right)\Gamma\left(\frac{s}{2}\right)}
\left(\frac{2t^{\frac{\beta}{\alpha}}}{|\mathrm{x}|}\right)^{-s}ds\, ,
\end{equation}
where $\frac{n}{2} - \frac{1}{2}<\gamma<\min (\alpha, n)$. Let us note that the Mellin-Barnes integral \eqref{eq20} can be interpreted as a particular case of the Fox H-function, too. The theory of the H-function, its properties, and applications were presented in a number of textbooks and papers (see e.g. \cite{Fox61}, \cite{Kiryakova}, \cite{Mai-P}, \cite{MatSax78}, \cite{YakLuc94}) so that here we do not discuss this subject  in detail and prefer to directly deduce the properties of the fundamental solution $G_{\alpha,\beta,n}$ from its Mellin-Barnes representation \eqref{eq20}. Starting with this representation and using simple linear variables substitutions, we can easily derive some other forms of this representation that will be useful for further discussions. Say, the substitutions $s \to -s $ and then $s \to s-n$ in the Mellin-Barnes representation \eqref{eq20} result in two other equivalent representations
\begin{equation}
\label{eq20_1}
G_{\alpha,\beta,n}(\mathrm{x},t)=\frac{1}{\alpha}\frac{|\mathrm{x}|^{-n}}{\pi^{\frac{n}{2}}}\frac{1}{2\pi i} \int_{\gamma-i\infty}^{\gamma+i\infty} \frac{\Gamma\left(\frac{n}{2}+\frac{s}{2}\right)\Gamma\left(-\frac{s}{\alpha}\right) \Gamma\left(1+\frac{s}{\alpha}\right)}{\Gamma\left(1+\frac{\beta}{\alpha}s\right)\Gamma\left(-\frac{s}{2}\right)}
\left(\frac{|\mathrm{x}|}{2t^{\frac{\beta}{\alpha}}}\right)^{-s}ds\, 
\end{equation}
and 
\begin{equation}
\label{eq20_2}
G_{\alpha,\beta,n}(\mathrm{x},t)=\frac{1}{\alpha}\frac{t^{-\frac{\beta n}{\alpha}}}{(4\pi)^{\frac{n}{2}}}\frac{1}{2\pi i} \int_{\gamma-i\infty}^{\gamma+i\infty} \frac{\Gamma\left(\frac{s}{2}\right)\Gamma\left(\frac{n}{\alpha}-\frac{s}{\alpha} \right) \Gamma\left(1-\frac{n}{\alpha}+\frac{s}{\alpha}\right)}
{\Gamma\left(1-\frac{\beta}{\alpha}n+\frac{\beta}{\alpha}s\right)\Gamma\left(\frac{n}{2}-\frac{s}{2}\right)}
\left(\frac{|\mathrm{x}|}{2t^{\frac{\beta}{\alpha}}}\right)^{-s}ds\, 
\end{equation}
that are valid under the conditions $-\min (\alpha, n)<\gamma<\frac{1}{2}-\frac{n}{2}$ and $\max (n-\alpha, 0)<\gamma<n$, respectively.  

It is worth mentioning that the Mellin-Barnes integrals at the right-hand sides of the  representations \eqref{eq20}, \eqref{eq20_1}, and \eqref{eq20_2} are well defined for $0<\alpha$, $0<\beta \le 2$, $n\in \mathbb{N}$ and thus the fundamental solution $G_{\alpha,\beta,n}(\mathrm{x},t)$ can be represented by these Mellin-Barnes integrals (at least) for $0<\alpha \le 2$, $0<\beta \le 2$, $n\in \mathbb{N}$. 

Finally, let us demonstrate how these integral representations can be used, say, for deriving some series representations of $G_{\alpha,\beta,n}(\mathrm{x},t)$ and then its representations in terms of elementary or special functions of the hypergeometric type. To this end, we consider a simple example. In the case $\beta = 1$ and $\alpha =2$ (standard diffusion equation), the representation \eqref{eq20_2}
takes the following form (two pairs of the Gamma-functions in the integral at the right-hand side of \eqref{eq20_2} are canceled):
$$
G_{2,1,n}(\mathrm{x},t)= \frac{t^{-\frac{n}{2}}}{2\, (4\pi)^{\frac{n}{2}}} \frac{1}{2\pi i }
\int_{\gamma-i\infty}^{\gamma+i\infty} \Gamma\left(\frac{ s}{2}\right)\, \left(\frac{z}{2}\right)^{-s}ds,\ z = \frac{|x|}{\sqrt{t}}. 
$$
Substitution of the variables $s \to 2s$ leads to an even simpler representation
$$
G_{2,1,n}(\mathrm{x},t)= \frac{t^{-\frac{n}{2}}}{(4\pi)^{\frac{n}{2}}} \frac{1}{2\pi i }
\int_{\gamma-i\infty}^{\gamma+i\infty} \Gamma\left(s\right)\, \left(\frac{z}{2}\right)^{-2s}ds,\ z = \frac{|x|}{\sqrt{t}}. 
$$
According to the Cauchy theorem, the contour of
integration in the integral at the right-hand side of the last formula can be transformed to the
loop $L_{-\infty}$ starting and ending at $-\infty$ and encircling
all poles $s_k=-k,\ k=0,1,2,\dots$ of the function
$\Gamma(s)$. Taking into account the Jordan lemma, the formula  for the residuals of the Gamma function (see e.g. \cite{Mar83})
$$
\mbox{res}_{s=-k} \Gamma(s) = \frac{(-1)^k}{k!},\ k=0,1,2,\dots
$$
and the Cauchy residue theorem lead to a series representation of 
$G_{2,1,n}(\mathrm{x},t)$:
$$
G_{2,1,n}(\mathrm{x},t)=\frac{t^{-\frac{n}{2}}}{(4\pi)^{\frac{n}{2}}}\int_{\gamma-i\infty}^{\gamma+i\infty} \Gamma(s)
\left(\frac{z}{2}\right)^{-2s}ds=\frac{t^{-\frac{n}{2}}}{(4\pi)^{\frac{n}{2}}}  \sum_{k=0}^\infty \frac{(-1)^k}{k!}  \left(\frac{z}{2}\right)^{2k}, \ z = \frac{|x|}{\sqrt{t}}.
$$
Thus the fundamental solution $G_{2,1,n}$ to the $n$-dimensional diffusion equation takes its well-known form:
\begin{equation}\label{eq23}
G_{2,1,n}(\mathrm{x},t)= 
\frac{1}{(\sqrt{4\pi t})^{n}}\exp\left(-\frac{|\mathrm{x}|^{2}}{4t}\right).
\end{equation} 

\section{Completely monotone functions and pdfs}

A very essential feature of the subordination formulas of  type \eqref{FS-S1} or \eqref{FS-S2} is that  their kernel functions $\Phi=  \Phi(s,t)$ can be interpreted as pdfs in $s,\ s\in \mathbb{R}_+$ for each value of $t,\ t>0$, i.e., that for any $t,\, s>0$ 
\begin{equation}
\label{3.1}
\Phi(s,t)\ge 0\ \mbox{ and } \ \int_0^\infty \Phi(s,t)\, ds = 1.
\end{equation} 

Verifying the properties \eqref{3.1} for a given special function $\Phi$ is often a very difficult task. In this section, a simple but efficient procedure will be suggested that helps to check  \eqref{3.1} for some special functions given in terms of the Mittag-Leffler function, the Wright function, and the generalized Wright function. We shell need these functions as kernels for the subordination formulas in the next section. This procedure uses the well-known connection between the non-negative functions and the complete monotone functions, but in the form written in terms of the Mellin integral transform. 

To start with, let us first give a definition of the completely monotone functions: 
\begin{Def}
A  non-negative function $ \phi:(0,\infty) \to \mathbb{R}$  is called a  completely
monotone function if it is of class $C^{\infty}(0,\infty)$  and
$(-1)^n \phi^{(n)}(\lambda)\geq 0$ for all $n \in \mathbb{N}$ and
$\lambda >0$.
\end{Def}




The functions $e^{-a\lambda^\alpha},\ a\ge 0,\ \alpha \le 1$ and
$E_{\alpha,\beta} (-\lambda),\ 0<\alpha \le 1,\ \alpha \le \beta$
are  well-known examples of completely monotone functions. Here
$E_{\alpha,\beta}(z)$ denotes the generalized Mittag-Leffler
function defined by the following convergent series
\begin{equation}
\label{Mit-Lef}
E_{\alpha,\beta}(z)=\sum_{k=0}^\infty \frac{z^k}{\Gamma(\alpha k
+ \beta)},\ \alpha >0, \ \beta \in \mathbb{C}.
\end{equation}

For more examples, properties, and
applications of the completely monotone functions we refer e.g., to  \cite{feller}, \cite{Miller-Samko}, and \cite{[SSV]}.

The basic property of the completely monotone
functions that we need in this section is the following one (Bernstein theorem): A
function $ \varphi:(0,\infty) \to \mathbb{R}$  is  completely monotone if and
only if it can be represented as the Laplace transform of a
non-negative measure (non-negative function or generalized
function). Because the Laplace transform is a Mellin convolution
type integral transform, the technique of the Mellin transform can be applied
for investigation of completely monotone functions.

Let the representation
\begin{equation}
\label{cm-1}
\phi(\lambda) = \int_0^\infty e^{-\lambda t}\Phi(t)\, dt,\ \lambda>0
\end{equation}
hold true for a non-negative function $\Phi$ with a known Mellin transform. Then the function $\phi$ is completely monotone and its Mellin transform is given by the
 formula (see e.g. \cite{LK2013})
\begin{equation}
\label{cm-2}
\phi^*(s) = \Gamma(s)\Phi^*(1-s)
\end{equation}
that can be transformed to the form
\begin{equation}
\label{cm-2-F} 
\Phi^*(s) = \frac{\phi^*(1-s)}{\Gamma(1-s)}\,.
\end{equation}
If the function $\Phi(t),\ t>0$ is non-negative, then the function
$\Phi_{\gamma,\beta}(t)= t^\gamma \Phi(t^{-\beta})$ is non-negative for any
$\gamma,\beta \in \mathbb{R}$, too. Thus the function $\phi_{\gamma,\beta}$ of the form
\begin{equation}
\label{cm-1-g}
\phi_{\gamma,\beta}(\lambda) = \int_0^\infty e^{-\lambda t}\Phi_{\gamma,\beta}(t)\, dt,\ \lambda>0
\end{equation}
is completely monotone and it follows from the relation (\ref{cm-2}) that
\begin{equation}
\label{cm-3}
\phi_{\gamma,\beta}^*(s) = \Gamma(s)\Phi_{\gamma,\beta}^*(1-s).
\end{equation}
Using the operational rules \eqref{(1.52)}-\eqref{(1.53)} for the Mellin integral transform, the Mellin transform of $\Phi_{\gamma,\beta}(t)= t^\gamma \Phi(t^{-\beta})$
can be written in the form
$$
\Phi_{\gamma,\beta}^*(s)= \frac{1}{|\beta|}\Phi^*\left( -\frac{\gamma}{\beta} - \frac{s}{\beta} \right).
$$
Thus  we get the following formula for $\phi_{\gamma,\beta}^*(s)$ defined by (\ref{cm-3}):
$$
\phi_{\gamma,\beta}^*(s) = \frac{1}{|\beta|} \Gamma(s) \Phi^*\left(\frac{s}{\beta} - \frac{1+\gamma}{\beta}   \right).
$$
The completely monotone function $\phi_{\gamma,\beta}$ given by (\ref{cm-1-g}) can be then represented
  as the Mellin-Barnes integral (inverse Mellin integral transform of $\phi_{\gamma,\beta}^*(s)$)
\begin{equation}
\label{cm-4} 
\phi_{\gamma,\beta}(\lambda) = \frac{1}{2\pi
i}\int_{\gamma-i\infty}^{\gamma+i\infty} \frac{1}{|\beta|} \Gamma(s)
\Phi^*\left( \frac{s}{\beta} -\frac{1+\gamma}{\beta} \right) \lambda^{-s}\,
ds.
\end{equation}
In many cases the function $\phi$ (and thus the function $\Phi$) is a particular case of the
Fox $H$-function and then $\Phi^*$ is represented in form of a quotient
of products of the Gamma functions. This means that the new completely
monotone function $\phi_{\gamma,\beta}$ given by \eqref{cm-4} is a particular case of the $H$-function, too.

Let us consider a simple example. It is known that the
exponential function $\phi(\lambda) = \exp(-\lambda^\alpha),\ 0<\alpha \le 1$ is
completely monotone. Its Mellin integral transform is given by the formula (\cite{LK2013},\ \cite{Mar83})
\begin{equation}
\label{(1.57)}
e^{-\lambda^\alpha}  \overset{\mathcal M}{\longleftrightarrow}  \frac{1}{\vert \alpha \vert}\Gamma(s/\alpha), \ 
 \Re(s/\alpha)>0.
\end{equation}
 The function $\Phi^*$ from (\ref{cm-2-F}) has then the
form
$$
\Phi^*(s) = \frac{\phi^*(1-s)}{\Gamma(1-s)} = \frac{1}{\alpha}
\frac{\Gamma\left( \frac{1}{\alpha}- \frac{s}{\alpha}
\right)}{\Gamma(1-s)}\,.
$$
It follows from the arguments presented above that the function
\begin{equation}
\label{cm-5}
\phi_{\gamma,\beta}(\lambda) = \frac{1}{2\pi i}\int_{\gamma-i\infty}^{\gamma+i\infty}
\frac{1}{\alpha |\beta|} \frac{\Gamma(s) \Gamma\left( \frac{\beta+\gamma+1}
{\alpha\beta}- \frac{s}{\alpha\beta} \right)}
{\Gamma\left( \frac{\beta+\gamma+1}{\beta}- \frac{s}{\beta} \right) } \lambda^{-s}\, ds
\end{equation}
is completely monotone, too.
 The function $\phi_{\gamma,\beta}$ given by (\ref{cm-5}) is evidently a particular case of the Fox $H$-function.
 In particular, in the case $\beta > \frac{1}{\alpha} -1$ it can be represented
 as the convergent series (see \cite{LK2013} or \cite{Mar83}) 
\begin{equation}
\label{cm-6} 
\phi_{\gamma,\beta}(\lambda) = \frac{1}{\alpha |\beta|} \sum_{k=0}^\infty
\frac{\Gamma\left(
\frac{\beta+\gamma+1}{\alpha\beta}+\frac{k}{\alpha\beta} \right)}
{k!\, \Gamma\left( \frac{\beta+\gamma+1}{\beta}+ \frac{k}{\beta}
\right) } (-\lambda)^{k}.
\end{equation}
We can easily recognize the last series as a particular case of the generalized Wright function
   defined by the series
\begin{equation}
\label{genWr}
\mbox{ }_p\Psi_q\left[{(a_1,A_1),\dots,(a_p,A_p)\atop (b_1,B_1)\dots
(b_q,B_q)};z\right]=\sum_{k=0}^\infty \frac{\prod_{i=1}^p \Gamma(a_i +A_i k)
}{\prod_{i=1}^q \Gamma(b_i+B_i k)} \frac{z^k}{ k!}
\end{equation}
for the $z$-values where 
  the series converges, and by the analytic continuation
of this series for other $z$-values. 
Thus, we have proved that the generalized Wright function
\begin{equation}
\label{cm-7}
\phi_{\gamma,\beta}(\lambda) = \frac{1}{\alpha |\beta|} \mbox{ }_1\Psi_1\left[{\left(\frac{\beta+\gamma+1}{\alpha\beta},\frac{1}
{\alpha\beta}\right)\atop \left(\frac{\beta+\gamma+1}{\beta},\frac{1}{\beta}\right)};-\lambda \right]
\end{equation}
is completely monotone under the conditions $0<\alpha \le 1,\ \frac{1}{\alpha} -1 < \beta$.
 In particular, let us set the following parameter values: $\beta = \frac{1}{\alpha},\ \gamma=-\frac{1}{\alpha}$.
 Then the series (\ref{cm-6}) (and thus the function \eqref{cm-7}) takes the form
\begin{equation}
\label{cm-8}
\phi_{\gamma,\beta}(\lambda) = \sum_{k=0}^\infty  \frac{\Gamma( 1+k)}
{k!\, \Gamma( \alpha+ \alpha k) } (-\lambda)^{k} =  \sum_{k=0}^\infty  \frac{(-\lambda)^{k}}
{\Gamma( \alpha+ \alpha k) }
\end{equation}
that defines the Mittag-Leffler function $E_{\alpha,\alpha}(-\lambda)$,
known to be completely monotone for $0<\alpha \le 1$. Taking
other known completely monotone functions and applying the procedure
described above, other new completely monotone functions can be
easily derived.

Another simple but important observation from the discussions presented above is that the Mellin integral transforms of the non-negative and completely monotone functions are connected by the formulas \eqref{cm-2} and \eqref{cm-2-F}. Say, if a function $\phi$ is completely monotone then the function $\Phi$ with the Mellin integral transform given by the
 formula 
$$
\Phi^*(s) = \frac{\phi^*(1-s)}{\Gamma(1-s)}
$$
is non-negative. Vise versa, if a function $\Phi$ is non-negative then the function $\phi$ with the Mellin integral transform given by the formula
$$
\phi^*(s) = \Gamma(s)\Phi^*(1-s)
$$
is completely monotone. 

Let us illustrate this procedure by some examples. 

\begin{Ex}
We start with a well-known pair of functions, namely, with the generalized Mittag-Leffler function defined by \eqref{Mit-Lef} and the Wright function defined by the convergent series
\begin{equation}
\label{Wr}
W_{a,\mu}(z)=\sum_{k=0}^\infty \frac{z^k}{k! \Gamma(a +\mu k)},\ \mu > -1,\ 
a,\, z\in \mathbb{C}
\end{equation}
and show that the function $p_{\alpha,\beta}(t) = \Gamma(\beta)\, W_{\beta-\alpha,-\alpha}(-t)$ can be interpreted as a pdf if $0<\alpha < 1,\ \alpha \le \beta$. 
\end{Ex}

The Mellin integral transforms of the generalized Mittag-Leffler function and of the Wright function are well-known (see e.g. \cite{LK2013}):
\begin{equation}
 \label{(1.67)}
E_{\alpha,\beta}(-t) \overset{\mathcal M}{\longleftrightarrow} 
\frac{\Gamma(s)\Gamma(1-s)}
{\Gamma(\beta -\alpha s)}
\end{equation}
$$
  \mbox{if}
 \ 0 < \Re(s) < 1,\ 0<\alpha < 2 \ \mbox{or}\   0 < \Re(s) < \min\{1,\Re(\beta)/2\}, \ \alpha=2, 
$$
\begin{equation}
\label{(1.67-wright)}
W_{a,\mu}(-t)  \overset{\mathcal M}{\longleftrightarrow} 
\frac{\Gamma(s)}
{\Gamma(a -\mu s)}
\end{equation}
$$
 \mbox{if}
 \ 0 < \Re(s),\ \mu < 1 \ \mbox{or}  \ 
0 < \Re(s) < \Re(a)/2 -1/4, \ \mu=1.
$$
As already mentioned, the function $\phi(\lambda) = E_{\alpha,\beta} (-\lambda)$ is completely monotone provided the conditions $0<\alpha \le 1,\ \alpha \le \beta$ are fulfilled. The Mellin integral transform of the function $\phi$ is given by 
\eqref{(1.67)}. Then the function $\Phi$ with the Mellin integral transform
$$
\Phi^*(s) = \frac{\phi^*(1-s)}{\Gamma(1-s)} = \frac{\Gamma(1-s)\Gamma(s)}
{\Gamma(1-s) \Gamma(\beta -\alpha + \alpha s)} = \frac{\Gamma(s)}
{\Gamma(\beta -\alpha + \alpha s)}
$$
is non-negative. Comparing this formula with \eqref{(1.67-wright)}, we conclude that the Wright function $W_{\beta-\alpha,-\alpha}(-t)$ is non-negative under the conditions $0<\alpha < 1,\ \alpha \le \beta$, i.e.,
\begin{equation}
 \label{Wr-non}
W_{\beta-\alpha,-\alpha}(-t)\ge 0,\ t>0,\ 0<\alpha < 1,\ \alpha \le \beta.
\end{equation}
Let us now check that the function $p_{\alpha,\beta}(t) = \Gamma(\beta)\, W_{\beta-\alpha,-\alpha}(-t)$ is a pdf on $\mathbb{R}_+$. Indeed, it is non-negative because of \eqref{Wr-non}. To calculate the integral of $p_{\alpha,\beta}$ over $\mathbb{R}_+$ let us mention that it can be interpreted as the Mellin integral transform of $p_{\alpha,\beta}$ at the point $s=1$. The formula \eqref{(1.67-wright)} leads now to the following chain of equalities:
$$
\int_0^\infty p_{\alpha,\beta}(t) \, dt = \int_0^\infty  \Gamma(\beta) W_{\beta-\alpha,-\alpha}(-t)\, dt =  \left. \frac{\Gamma(\beta) \Gamma(s)}
{\Gamma(\beta-\alpha + \alpha s)}\right|_{s=1} = \frac{\Gamma(\beta) }
{\Gamma(\beta)} = 1.
$$

\begin{Ex}
In this example, we verify that the following function defined in terms of the Mellin-Barnes integral
\begin{equation}
 \label{Ex3_1}
 \Phi_{\alpha,\beta}(t) = \frac{2}{\alpha} \frac{1}{2\pi i}\int_{\gamma-i\infty}^{\gamma+i\infty} \frac{\Gamma\left( \frac{2}{\alpha}- \frac{2}{\alpha}s\right) \Gamma\left( 1-\frac{2}{\alpha}+ \frac{2}{\alpha}s\right)}
 {\Gamma\left( 1-\frac{2\beta}{\alpha} + \frac{2\beta}{\alpha}s\right) \Gamma\left( 1-s\right)}\, t^{-s}\, ds
 \end{equation}
 can be interpreted as a pdf on $\mathbb{R}_+$ for $0 < \beta \le 1$ and $0<\alpha \le 2$ when $\alpha +2\beta <4$.  
\end{Ex}
According the the general theory of the Mellin-Barnes integrals (see e.g. \cite{Mar83}), the Mellin-Barnes integral \eqref{Ex3_1} exists for $\frac{2}{\alpha} -1 <\Re(s) < \frac{2}{\alpha}$ under the conditions $0 < \beta$, $0<\alpha$,  and $\alpha +2\beta <4$ and its Mellin transform can be calculated as follows:
\begin{equation}
 \label{Ex3_2}
 \Phi_{\alpha,\beta}^*(s)  = \frac{2}{\alpha} \frac{\Gamma\left( \frac{2}{\alpha}- \frac{2}{\alpha}s\right) \Gamma\left( 1-\frac{2}{\alpha}+ \frac{2}{\alpha}s\right)}
 {\Gamma\left( 1-\frac{2\beta}{\alpha} + \frac{2\beta}{\alpha}s\right) \Gamma\left( 1-s\right)}.
 \end{equation}
Now we construct the function $\phi^*(s)$ according to the formula \eqref{cm-2}:
\begin{equation}
 \label{Ex3_27}
\phi^*(s) = \Gamma(s)\Phi_{\alpha,\beta}^*(1-s) = \Gamma(s) \frac{2}{\alpha} \frac{\Gamma\left(\frac{2}{\alpha}s\right) \Gamma\left( 1-\frac{2}{\alpha}s\right)}
 {\Gamma\left( 1-\frac{2\beta}{\alpha}s\right) \Gamma\left(s\right)} =
\frac{2}{\alpha} \frac{\Gamma\left(\frac{2}{\alpha}s\right) \Gamma\left( 1-\frac{2}{\alpha}s\right)}
 {\Gamma\left( 1-\frac{2\beta}{\alpha}s\right)}. 
\end{equation}
The function $\phi=\phi(\lambda)$ can be then represented as the following Mellin-Barnes integral:
$$
\phi(\lambda) =  \frac{1}{2\pi i}\int_{\gamma-i\infty}^{\gamma+i\infty} \frac{2}{\alpha} \frac{\Gamma\left(\frac{2}{\alpha}s\right) \Gamma\left( 1-\frac{2}{\alpha}s\right)}
 {\Gamma\left( 1-\frac{2\beta}{\alpha}s\right)}\, \lambda^{-s}\, ds.
 $$
The variables substitution $\frac{2}{\alpha}s \to s$ leads to the representation
$$
\phi(\lambda) =  \frac{1}{2\pi i}\int_{\gamma-i\infty}^{\gamma+i\infty} \frac{\Gamma\left(s\right) \Gamma\left( 1- s\right)}
 {\Gamma\left( 1- \beta s\right)}\, \left(\lambda^{\frac{\alpha}{2}}\right)^{-s}\, ds.
 $$
Comparing this formula with the Mellin transform \eqref{ML-Mellin1} of the Mittag-Leffler function, we arrive at the representation
\begin{equation}
 \label{Ex3_3}
\phi(\lambda) = E_\beta \left(-\lambda^{\frac{\alpha}{2}}\right).
\end{equation}
The Mittag-Leffler function $f(\lambda)=E_\beta(-\lambda)$ is known to be completely monotone for $0 <\beta \le 1$. Thus for $\alpha =2$ the function $\phi(\lambda)$ defined by \eqref{Ex3_3} is completely monotone. Now let $\alpha$ satisfy the inequalities $0<\alpha <2$. Then the function $g(\lambda)=\lambda^{\frac{\alpha}{2}}$ is a Bernstein function because its derivative $g^\prime(\lambda)=\frac{\alpha}{2}\lambda^{\frac{\alpha}{2}-1}$ is completely monotone. But a composition of a completely monotone function and a Bernstein function is completely monotone (see e.g. [SSV]). Thus the function $\phi(\lambda) = f(g(\lambda))$ is completely monotone for $0<\alpha <2$, too. Because 
$\phi^*(s)$ and $\Phi_{\alpha,\beta}^*(s)$ are connected by the formula \eqref{cm-2} and the function $\phi$ is completely monotone, it follows now that $\Phi_{\alpha,\beta}(t)$ is non-negative, i.e.,
$$
\Phi_{\alpha,\beta}(t) \ge 0,\ t>0,\ 0 < \beta \le 1, \ 0<\alpha \le 2,\ \alpha +2\beta <4.
$$
To evaluate the integral of $\Phi_{\alpha,\beta}(t)$ over $\mathbb{R}_+$ we again use the technique of the Mellin integral transform:
$$
\int_0^\infty \Phi_{\alpha,\beta}(t) \, dt =  \lim_{s\to 1}  \frac{2}{\alpha} \frac{\Gamma\left( \frac{2}{\alpha}(1-s)\right) \Gamma\left( 1-\frac{2}{\alpha}+ \frac{2}{\alpha}s\right)}
 {\Gamma\left( 1-\frac{\beta}{\alpha} + \frac{2\beta}{\alpha}s\right) \Gamma\left( 1-s\right)} = \frac{2}{a} \lim_{s\to 1} \frac{\Gamma\left( \frac{2}{\alpha}(1-s)\right)}
 { \Gamma\left( 1-s\right)}   =1.
$$

\section{Subordination formulas for the fundamental solution}

To demonstrate our method, we open this section with derivation of some known subordination formulas of type \eqref{1.3} for the fundamental solution $G_{\alpha,\beta, n}$. The starting point is the Mellin-Barnes representation \eqref{eq20_2} that we rewrite in the form 
\begin{equation}
\label{MBR}
G_{\alpha,\beta,n}(\mathrm{x},t)=\frac{1}{\alpha}\frac{t^{-\frac{\beta n}{\alpha}}}{(4\pi)^{\frac{n}{2}}}\frac{1}{2\pi i} \int_{\gamma-i\infty}^{\gamma+i\infty} 
K_{\alpha,\beta,n}(s) 
z^{-s}ds\, , \ \ z = \frac{|\mathrm{x}|}{2t^{\frac{\beta}{\alpha}}}
\end{equation}
with 
\begin{equation}
\label{K}
K_{\alpha,\beta,n}(s) = \frac{\Gamma\left(\frac{s}{2}\right)\Gamma\left(\frac{n}{\alpha}-\frac{s}{\alpha} \right) \Gamma\left(1-\frac{n}{\alpha}+\frac{s}{\alpha}\right)}
{\Gamma\left(1-\frac{\beta}{\alpha}n+\frac{\beta}{\alpha}s\right)\Gamma\left(\frac{n}{2}-\frac{s}{2}\right)}.
\end{equation}
By setting $\beta = 1$ in the formulas above we get the fundamental solution of the space-fractional diffusion equation in the form:
\begin{equation}
\label{MBR_1}
G_{\alpha,1,n}(\mathrm{x},t)=\frac{1}{\alpha}\frac{t^{-\frac{n}{\alpha}}}{(4\pi)^{\frac{n}{2}}}\frac{1}{2\pi i} \int_{\gamma-i\infty}^{\gamma+i\infty} 
K_{\alpha,1,n}(s) 
z^{-s}ds\, , \ \ z = \frac{|\mathrm{x}|}{2t^{\frac{1}{\alpha}}}
\end{equation}
with 
\begin{equation}
\label{K_1}
K_{\alpha,1,n}(s) = \frac{\Gamma\left(\frac{s}{2}\right)\Gamma\left(\frac{n}{\alpha}-\frac{s}{\alpha} \right) }
{\Gamma\left(\frac{n}{2}-\frac{s}{2}\right)}.
\end{equation}
The key point for derivation of a subordination formula for $G_{\alpha,\beta,n}$ with $0<\beta <1$ is in observation that the kernel function $K_{\alpha,\beta,n}$ in the Mellin-Barnes integral \eqref{MBR} can be represented as product of two factors:
\begin{equation}
\label{Pr_1}
K_{\alpha,\beta,n}(s) = K_{\alpha,1,n}(s) \times \Phi_{\alpha,\beta,n}^*(s),
\end{equation}
where $K_{\alpha,1,n}(s)$ is the kernel function in the Mellin-Barnes integral \eqref{MBR_1} for the fundamental solution $G_{\alpha,1,n}$ and 
\begin{equation}
\label{Phi_1}
\Phi_{\alpha,\beta,n}^*(s) = 
\frac{\Gamma\left(1-\frac{n}{\alpha}+\frac{s}{\alpha}\right)} {\Gamma\left(1-\frac{\beta}{\alpha}n+\frac{\beta}{\alpha}s\right)}.
\end{equation}
Due to the convolution formula \eqref{eq17} for the Mellin transform, the product formula \eqref{Pr_1} in the Mellin domain leads to an integral representation of $G_{\alpha,\beta,n}$ in the form
\begin{equation}
\label{sub_1} 
G_{\alpha,\beta,n}(\mathrm{x},t) = \frac{1}{\alpha}\frac{t^{-\frac{\beta n}{\alpha}}}{(4\pi)^{\frac{n}{2}}} \int_0^\infty \Phi_{\alpha,\beta,n}(\tau) \tilde{G}_{\alpha,1,n}\left(\frac{z}{\tau}\right)\, \frac{d\tau}{\tau}, \ z = \frac{|\mathrm{x}|}{2t^{\frac{\beta}{\alpha}}},
\end{equation}
where $\Phi_{\alpha,\beta,n}(\tau)$ is the inverse Mellin integral transform of $\Phi_{\alpha,\beta,n}^*(s)$ given by \eqref{Phi_1} and 
\begin{equation}
\label{tilde_1} 
\tilde{G}_{\alpha,1,n}(\tau) = \frac{1}{2\pi i} 
\int_{\gamma-i\infty}^{\gamma+i\infty} 
K_{\alpha,1,n}(s) \, \tau^{-s}ds
\end{equation}
is a slightly modified fundamental solution $G_{\alpha,1,n}$:
\begin{equation}
\label{G_G_1}
G_{\alpha,1,n}(\mathrm{x},t)=\frac{1}{\alpha}\frac{t^{-\frac{n}{\alpha}}}{(4\pi)^{\frac{n}{2}}} \tilde{G}_{\alpha,1,n}(z), \ \ z = \frac{|\mathrm{x}|}{2t^{\frac{1}{\alpha}}}.
\end{equation}
The formula \eqref{sub_1} is  a subordination formula for the fundamental solution $G_{\alpha,\beta,n}$ and now we put it into the standard form. To do this, let us derive an explicit representation for the kernel function $\Phi_{\alpha,\beta,n}(\tau)$ that  is defined  as the  Mellin-Barnes integral (under the condition $0 <\beta <1$)
\begin{equation}
\label{Phi_MB_1}
\Phi_{\alpha,\beta,n}(\tau) = \frac{1}{2\pi i} 
\int_{\gamma-i\infty}^{\gamma+i\infty} 
\Phi_{\alpha,\beta,n}^*(s) \, \tau^{-s}\, ds = \frac{1}{2\pi i} 
\int_{\gamma-i\infty}^{\gamma+i\infty} \frac{\Gamma\left(1-\frac{n}{\alpha}+\frac{s}{\alpha}\right)} {\Gamma\left(1-\frac{\beta}{\alpha}n+\frac{\beta}{\alpha}s\right)}
\, \tau^{-s}\, ds.
\end{equation}
The general theory of the Mellin-Barnes integrals (see e.g. \cite{Mar83}) says that the contour of
integration in the integral at the right-hand side of \eqref{MBR_1} can be transformed to the
loop $L_{-\infty}$ starting and ending at $-\infty$ and encircling
all poles of the function
$\Gamma\left(1-\frac{n}{\alpha}+\frac{s}{\alpha}\right)$. Taking into account the Jordan lemma and  the formula for the residuals of the Gamma-function, 
 the Cauchy residue theorem leads to the following series representation of $\Phi_{\alpha,\beta,n}$ (for details we refer the reader to \cite{BL1} or \cite{L2017}):
\begin{equation}
\label{Phi_ser_1}
\Phi_{\alpha,\beta,n}(\tau) = \sum_{k=0}^\infty \frac{\alpha (-1)^k}{k!}
\frac{1}{\Gamma\left( 1-\beta -\beta k \right)} (\tau^\alpha)^{k+1-\frac{n}{\alpha}}
\end{equation}
that can be recognized to be a special case of the Wright function:
\begin{equation}
\label{Phi_wr_1}
\Phi_{\alpha,\beta,n}(\tau) = \alpha \tau^{\alpha- n} W_{1-\beta,-\beta}(-\tau^\alpha),\ 0 <\beta <1. 
\end{equation}
Putting now the formulas \eqref{G_G_1} and \eqref{Phi_wr_1} into the integral representation \eqref{sub_1} and substituting the variables $\tau^\alpha \to \tau$, we first get the formula
\begin{equation}
\label{sub_1_1} 
G_{\alpha,\beta,n}(\mathrm{x},t) = \int_0^\infty W_{1-\beta,-\beta}(-\tau) G_{\alpha,1,n}(\mathrm{x},t^\beta\tau)\, d\tau, 
\end{equation}
that can be transformed into the known subordination formula (see \eqref{1.1} with $\delta = 1$)
\begin{equation}
\label{sub_1_1_1} 
G_{\alpha,\beta,n}(\mathrm{x},t) = \int_0^\infty t^{-\beta} W_{1-\beta,-\beta}(-st^{-\beta}) G_{\alpha,1,n}(\mathrm{x},s)\, ds, \ 0<\beta <1 
\end{equation}
by the variables substitution $t^\beta \tau \to s$. 

Let us note here that we can express the fundamental solution $G_{\alpha,1,n}$ in terms of the generalized Wright function $_1\Psi_1$. Indeed, applying the same technique as for the Mellin-Barnes integral \eqref{Phi_MB_1},
we first get a series representation of $G_{\alpha,1,n}$:
\begin{equation}
\label{Ser_1}
G_{\alpha,1,n}(\mathrm{x},t) = \frac{1}{\alpha}\frac{t^{-\frac{n}{\alpha}}}{(4\pi)^{\frac{n}{2}}} \sum_{k=0}^\infty \frac{2(-1)^k}{k!}
\frac{\Gamma\left( \frac{n}{\alpha}+ \frac{2}{\alpha}k \right)} {\Gamma\left( \frac{n}{2}+ k \right)} \left(\frac{|\mathrm{x}|^2}{4t^{\frac{2}{\alpha}}}\right)^k. 
\end{equation}
Comparing this series with \eqref{genWr}, we can represent $G_{\alpha,1,n}$  in  
terms of the generalized Wright function:
\begin{equation}
\label{Wri_1}
G_{\alpha,1,n}(\mathrm{x},t) =
\frac{2}{\alpha}\frac{t^{-\frac{n}{\alpha}}}{(4\pi)^{\frac{n}{2}}} 
\mbox{ }_1\Psi_1\left[{\left(\frac{n}{\alpha},\frac{2}
{\alpha}\right)\atop \left(\frac{n}{2},1\right)};
-\frac{|\mathrm{x}|^2}{4t^{\frac{2}{\alpha}}} \right]. 
\end{equation}
It is worth mentioning that the generalized Wright function from the right-hand side of the formula \eqref{Wri_1} is a particular case of the function \eqref{cm-7} and thus completely monotone with respect to the variable $z = \frac{|\mathrm{x}|^2}{4t^{\frac{2}{\alpha}}}$. Because the function $W_{1-\beta,-\beta}(-t)$ is non-negative if $0<\beta <1$ (see Example 1 of the previous section), the subordination formula \eqref{sub_1_1_1} along with the representation \eqref{Wri_1} means that the fundamental solution $G_{\alpha,\beta,n}$ is non-negative for $0<\beta <1$, $0<\alpha \le 2$ and it is a pdf in $\mathrm{x}$ for each $t>0$ that can be easily shown by the technique of the Mellin integral transform.

Now we consider the two-dimensional $\alpha$-fractional diffusion equation that is obtained from \eqref{eq1} for the parameter values $n=2$ and $\beta = \alpha/2$ (for derivation of this equation and analysis of its mathematical, physical, and probabilistic properties see \cite{L2016}). 

Specializing the formulas \eqref{MBR} and \eqref{K} for this case, we obtain the representations:
\begin{equation}
\label{MBR_2}
G_{\alpha,\alpha/2,2}(\mathrm{x},t)=\frac{1}{\alpha}\frac{t^{-1}}{4\pi}\frac{1}{2\pi i} \int_{\gamma-i\infty}^{\gamma+i\infty} 
K_{\alpha,\alpha/2,2}(s) 
z^{-s}ds\, , \ \ z = \frac{|\mathrm{x}|}{2t^{\frac{1}{2}}}
\end{equation}
with 
\begin{equation}
\label{K_2}
K_{\alpha,\alpha/2,2}(s) = \frac{\Gamma\left(\frac{2}{\alpha}-\frac{s}{\alpha}\right)\Gamma\left(1-\frac{2}{\alpha}+\frac{s}{\alpha}\right) }
{\Gamma\left(1-\frac{s}{2}\right)}.
\end{equation}
Now we consider the kernel function $K_{\alpha,\beta,2}(s)$ (the function \eqref{K} with $n=2$) under the condition $\beta < \frac{\alpha}{2}$ and represent it as follows:
\begin{equation}
\label{Pr_2}
K_{\alpha,\beta,2}(s) = K_{\alpha,\alpha/2,2}(s) \times \Phi_{\alpha,\beta}^*(s),
\end{equation}
where $K_{\alpha,\alpha/2,2}(s)$ is the kernel function in the Mellin-Barnes integral \eqref{MBR_2} for the fundamental solution $G_{\alpha,\alpha/2,2}$ and 
\begin{equation}
\label{Phi_2}
\Phi_{\alpha,\beta}^*(s) = 
\frac{\Gamma\left(\frac{s}{2}\right)} {\Gamma\left(1-\frac{2 \beta}{\alpha}+\frac{\beta}{\alpha}s\right)}.
\end{equation}
Once again, the product formula \eqref{Pr_2} in the Mellin domain  leads to an integral representation of $G_{\alpha,\beta,2}$ in the form
\begin{equation}
\label{sub_2} 
G_{\alpha,\beta,2}(\mathrm{x},t) = \frac{1}{\alpha}\frac{t^{-\frac{2 \beta }{\alpha}}}{4\pi} \int_0^\infty \Phi_{\alpha,\beta}(\tau) \tilde{G}_{\alpha,\alpha/2,2}\left(\frac{z}{\tau}\right)\, \frac{d\tau}{\tau}, \ z = \frac{|\mathrm{x}|}{2t^{\frac{\beta}{\alpha}}},
\end{equation}
where $\Phi_{\alpha,\beta}(\tau)$ is the inverse Mellin integral transform of $\Phi_{\alpha,\beta}^*(s)$ given by \eqref{Phi_2} and 
\begin{equation}
\label{tilde_2} 
\tilde{G}_{\alpha,\alpha/2,2}(\tau) = \frac{1}{2\pi i} 
\int_{\gamma-i\infty}^{\gamma+i\infty} 
K_{\alpha,\alpha/2,2}(s) \, \tau^{-s}ds
\end{equation}
is connected with  $G_{\alpha,\alpha/2,2}$ by the relation
\begin{equation}
\label{G_G_2}
G_{\alpha,\alpha/2,2}(\mathrm{x},t)=\frac{1}{\alpha}\frac{t^{-1}}{4\pi} \tilde{G}_{\alpha,\alpha/2,2}(\tilde{z}), \ \ \tilde{z} = \frac{|\mathrm{x}|}{2t^{\frac{1}{2}}}.
\end{equation}
The formula \eqref{sub_2} is  a subordination formula for the fundamental solution $G_{\alpha,\beta,2}$. To put it into the standard form, we first derive an explicit representation for the kernel function $\Phi_{\alpha,\beta}(\tau)$ that is defined as the  Mellin-Barnes integral
\begin{equation}
\label{Phi_MB_2}
\Phi_{\alpha,\beta}(\tau) = \frac{1}{2\pi i} 
\int_{\gamma-i\infty}^{\gamma+i\infty} 
\Phi_{\alpha,\beta}^*(s) \, \tau^{-s}\, ds = \frac{1}{2\pi i} 
\int_{\gamma-i\infty}^{\gamma+i\infty} \frac{\Gamma\left(\frac{s}{2}\right)} {\Gamma\left(1-\frac{2 \beta}{\alpha}+\frac{\beta}{\alpha}s\right)}
\, \tau^{-s}\, ds.
\end{equation}
Proceeding as above, for $\beta < \frac{\alpha}{2}$ we get (for details we refer the reader to \cite{L2016}):
\begin{equation}
\label{Phi_ser_2}
\Phi_{\alpha,\beta}(\tau) = \sum_{k=0}^\infty \frac{2(-1)^k}{k!}
\frac{1}{\Gamma\left( 1-\frac{2 \beta}{\alpha} -\frac{2 \beta}{\alpha} k \right)} (\tau^2)^{k}
\end{equation}
that can be recognized to be a special case of the Wright function:
\begin{equation}
\label{Phi_wr_2}
\Phi_{\alpha,\beta}(\tau) = 2 W_{1-\frac{2 \beta}{\alpha},-\frac{2 \beta}{\alpha}}(-\tau^2). 
\end{equation}
Putting now the formulas \eqref{G_G_2} and \eqref{Phi_wr_2} into the integral representation \eqref{sub_2} and after some elementary transformations, we get the subordination formula
\begin{equation}
\label{sub_1_1_2} 
G_{\alpha,\beta,2}(\mathrm{x},t) = \int_0^\infty t^{-\frac{2 \beta}{\alpha}} W_{1-\frac{2 \beta}{\alpha},-\frac{2 \beta}{\alpha}}(-st^{-\frac{2 \beta}{\alpha}}) G_{\alpha,\alpha/2,2}(\mathrm{x},s)\, ds,\ \beta <\frac{\alpha}{2}. 
\end{equation}
This formula is evidently a particular case of the subordination formula \eqref{1.1} with $n=2$ and $\delta = \alpha/2$. An advantage of our approach is that we can deduce a nice closed form formula for the fundamental solution $G_{\alpha,\alpha/2,2}$ being a part of  the formula \eqref{sub_1_1_2}. 

Again, we start with the Mellin-Barnes integral \eqref{MBR_2} with the kernel \eqref{K_2} and first obtain its series representation:
\begin{equation}
\label{Ser_2}
G_{\alpha,\alpha/2,2}(\mathrm{x},t) = \frac{1}{\alpha}\frac{t^{-1}}{4\pi} \sum_{k=0}^\infty \frac{\alpha(-1)^k}{k!}
\frac{\Gamma( 1+k )} {\Gamma\left( \frac{\alpha}{2}+ \frac{\alpha}{2} k \right)} \left(\frac{|\mathrm{x}|}{2t^{\frac{1}{2}}}\right)^{\alpha k +\alpha - 2} 
\end{equation}
that can be rewritten in terms of the generalized Mittag-Leffler function:
\begin{equation}
\label{Wri_2}
G_{\alpha,\alpha/2,2}(\mathrm{x},t) =
\frac{1}{4 \pi t}  \left(\frac{|x|}{2\sqrt{t}}\right)^{\alpha-2} \, E_{\frac{\alpha}{2},\frac{\alpha}{2}} \left(-\left(\frac{|x|}{2\sqrt{t}}\right)^\alpha\right). 
\end{equation}
The generalized Mittag-Leffler function from the right-hand side of the formula \eqref{Wri_2} is completely monotone with respect to the variable $z = \left(\frac{|x|}{2\sqrt{t}}\right)^\alpha$ and thus non-negative. Moreover, for each $t>0$ the fundamental solution $G_{\alpha,\alpha/2,2}(\mathrm{x},t)$ is a pdf in $\mathrm{x}$ (see \cite{L2016} for details).   

To obtain the subordination formula of type \eqref{1.1}, we compare the kernel functions $K_{\alpha,\beta,n}(s)$ and $K_{\alpha,\delta,n}(s)$ defined by the formula \eqref{K} with $0 <\beta <\delta \le 2$. Evidently, we can 
 represent $K_{\alpha,\beta,n}(s)$ as product of two factors:
\begin{equation}
\label{Pr_31}
K_{\alpha,\beta,n}(s) = K_{\alpha,\delta,n}(s) \times \Phi_{\alpha,\beta,n}^*(s),
\end{equation}
where 
\begin{equation}
\label{Phi_31}
\Phi_{\alpha,\beta,n}^*(s) = 
\frac{\Gamma\left(1-\frac{\delta n}{\alpha}+\frac{\delta}{\alpha}s\right)} {\Gamma\left(1-\frac{\beta n}{\alpha}+\frac{\beta}{\alpha}s\right)}.
\end{equation}
The function $\Phi_{\alpha,\beta,n}(\tau)$ can be determined as the inverse Mellin integral transform of $\Phi_{\alpha,\beta,n}^*(s)$  and then represented as a series 
\begin{equation}
\label{Phi_ser_31}
\Phi_{\alpha,\beta,n}(\tau) = \sum_{k=0}^\infty \frac{\alpha}{\delta} \frac{(-1)^k} 
{k!}
\frac{1}{\Gamma\left( 1-\frac{\beta}{\delta} -\frac{\beta}{\delta} k \right)} (\tau)^{\frac{\alpha}{\delta}k+\frac{\alpha}{\delta}-n}
\end{equation}
that can be recognized to be a special case of the Wright function:
\begin{equation}
\label{Phi_wr_31}
\Phi_{\alpha,\beta,n}(\tau) = \frac{\alpha}{\delta} \tau^{\frac{\alpha}{\delta}- n} W_{1-\frac{\beta}{\delta},-\frac{\beta}{\delta}}(-\tau^{\frac{\alpha}{\delta}}). 
\end{equation}
Repeating the argumentation that was employed for derivation of the subordination formulas \eqref{sub_1_1_1} and \eqref{sub_1_1_2}, we arrive at the subordination formula of type \eqref{1.1}. 

Now let us apply the method described above with respect to both the order $\alpha$ of the space-fractional derivative and the order $\beta$ of the time-fractional derivative. In our derivations, the four parameters Wright function in the form 
\begin{equation}
\label{wright}
W_{(a,\mu),(b,\nu)}(z)=\sum_{k=0}^\infty \frac{z^k}{\Gamma(a
+\mu k) \Gamma(b+\nu k)},\ \ 
\mu,\nu\in \mathbb{R},\ a,\, b,\, z\in \mathbb{C}
\end{equation}
will be used. This function was introduced in \cite{Wr2} for the positive values of the parameters $\mu$ and 
$\nu>0$. When $a=\mu=1$ or $b=\nu=1$, respectively, the four parameters Wright function 
is reduced to the Wright function
(\ref{Wr}). In \cite{Luc}, the four parameters Wright function was investigated in the case when one of the parameters $\mu$ or $\nu$ is negative.
In particular, it was proved there that the function $W_{(a,\mu),(b,\nu)}(z)$
is an entire function provided that $0<\mu+\nu,\ a,b\in \mathbb{C}$.


In the  case $\mu + \nu =0$, the four parameters Wright function is not en entire function anymore. The convergence radius of the series from \eqref{wright} with $\mu + \nu =0$ is equal to one, not to infinity, as can be seen from the asymptotics of the series terms as $k \to \infty$:
$$
\left| \frac{1}{ \Gamma(a
-\nu k) \Gamma(b+\nu k)}\right| = \left| \frac{\sin(\pi(a
-\nu k))}{\pi}\frac{\Gamma(1-a
+\nu k)}{ \Gamma(b+\nu k)}\right| =
$$
$$
 = \left| \frac{\cosh(\pi\Im(a))}{\pi} (\nu k)^{1-a-b} \left[ 1 + O(k^{-1})\right]\right|, \ k\to +\infty. 
$$

Now we formulate and prove the main result of this paper.

\begin{The}
\label{the1}
For the fundamental solution $G_{\alpha,\beta,n}(\mathrm{x},t)$ to the multi-dimensional space-time-fractional diffusion-wave equation \eqref{eq1} with $0 <\beta \le 1$, $0 < \alpha \le 2$, and $2\beta + \alpha <4$ the following subordination formula is valid:
\begin{equation}
\label{sub_1_1_3} 
G_{\alpha,\beta,n}(\mathrm{x},t) = \int_0^\infty t^{-\frac{2\beta}{\alpha}}\Phi_{\alpha,\beta}(st^{-\frac{2\beta}{\alpha}})\, G_{2,1,n}(\mathrm{x},s)\, ds, 
\end{equation}
where the fundamental solution to the conventional diffusion-wave equation 
is given by
$$
G_{2,1,n}(\mathrm{x},t) = \frac{1}{(\sqrt{4\pi t})^{n}}\exp\left(-\frac{|\mathrm{x}|^{2}}{4t}\right)
$$
and the kernel function $\Phi_{\alpha,\beta}(\tau)$ is a probability density function that is defined as follows:
\begin{equation}
\label{Phi_all} 
\Phi_{\alpha,\beta}(\tau)= 
\begin{cases}
\tau^{\frac{\alpha}{2}-1}\, W_{(1-\beta,-\beta),(\frac{\alpha}{2},\frac{\alpha}{2})}\left(-\tau^{\frac{\alpha}{2}}\right) \mbox{ if }\ \frac{\beta}{\alpha} <\frac{1}{2}, \\
\ \ \\
\tau^{-1}\, W_{(1,\beta),(0,-\frac{\alpha}{2})}\left(-\tau^{-\frac{\alpha}{2}}\right) \mbox{ if }\ \frac{\beta}{\alpha} > \frac{1}{2}, \\
\ \ \\
\begin{cases}
\frac{\tau^{\frac{\alpha}{2}-1}}{\pi} \sum_{k=0}^\infty \sin\left(  \frac{\pi\alpha}{2}(k+1) \right) \left(-\tau^{\frac{\alpha}{2}}\right)^k \mbox{ if }\ 0 < \tau < 1 \\
\ \ \\
-\frac{\tau^{-1}}{\pi} \sum_{k=0}^\infty \sin\left( \frac{\pi\alpha}{2}k \right) \left(-\tau^{-\frac{\alpha}{2}}\right)^k \mbox{ if }\ \tau > 1 \\
\end{cases}
\mbox{ if }\ \frac{\beta}{\alpha} = \frac{1}{2}.
\end{cases}
\end{equation}
 \end{The}

The method of derivation of the formula \eqref{sub_1_1_3} is the same as above. We start by putting $\alpha = 2$ and $\beta = 1$ into the the formulas \eqref{MBR} and \eqref{K} and  obtain a Mellin-Barnes representation for the fundamental solution to the conventional diffusion equation:
\begin{equation}
\label{MBR_3}
G_{2,1,n}(\mathrm{x},t)=\frac{1}{2}\frac{t^{-\frac{n}{2}}} 
{(4\pi)^{\frac{n}{2}}}
\frac{1}{2\pi i} \int_{\gamma-i\infty}^{\gamma+i\infty} 
K_{2,1,n}(s) 
z^{-s}ds\, , \ \ z = \frac{|\mathrm{x}|}{2t^{\frac{1}{2}}}
\end{equation}
with 
\begin{equation}
\label{K_3}
K_{2,1,n}(s) = \Gamma\left(\frac{s}{2}\right).
\end{equation}
The kernel function $K_{\alpha,\beta,n}(s)$ defined by \eqref{K}  can be then represented as follows:
\begin{equation}
\label{Pr_3}
K_{\alpha,\beta,n}(s) = K_{2,1,n}(s) \times \Psi_{\alpha,\beta,n}^*(s),
\end{equation}
where $K_{2,1,n}(s)$ is the kernel function in the Mellin-Barnes integral \eqref{MBR_3} for the fundamental solution $G_{2,1,n}$ and 
\begin{equation}
\label{Phi_3}
\Psi_{\alpha,\beta}^*(s) = 
\frac{\Gamma\left(\frac{n}{\alpha}-\frac{s}{\alpha} \right) \Gamma\left(1-\frac{n}{\alpha}+\frac{s}{\alpha}\right)}
{\Gamma\left(1-\frac{\beta}{\alpha}n+\frac{\beta}{\alpha}s\right)\Gamma\left(\frac{n}{2}-\frac{s}{2}\right)}.
\end{equation}
Because of the Mellin convolution theorem, the Mellin-Barnes integral \eqref{MBR} and the product formula \eqref{Pr_3} in the Mellin domain  lead to the integral representation 
\begin{equation}
\label{sub_3} 
G_{\alpha,\beta,n}(\mathrm{x},t) = \frac{1}{\alpha}\frac{t^{-\frac{\beta n}{\alpha}}}{(4\pi)^{\frac{n}{2}}} \int_0^\infty \Psi_{\alpha,\beta,n}(\tau) \tilde{G}_{2,1,n}\left(\frac{z}{\tau}\right)\, \frac{d\tau}{\tau}, \ z = \frac{|\mathrm{x}|}{2t^{\frac{\beta}{\alpha}}},
\end{equation}
where $\Psi_{\alpha,\beta,n}(\tau)$ is the inverse Mellin integral transform of $\Psi_{\alpha,\beta}^*(s)$ given by \eqref{Phi_3} and 
\begin{equation}
\label{tilde_3} 
\tilde{G}_{2,1,n}(\tau) = \frac{1}{2\pi i} 
\int_{\gamma-i\infty}^{\gamma+i\infty} 
K_{2,1,n}(s) \, \tau^{-s}ds.
\end{equation}
Comparing  \eqref{MBR_3} and \eqref{tilde_3}, we first get the relation
\begin{equation}
\label{G_G_31}
G_{2,1,n}(\mathrm{x},t)= \frac{1}{2}\frac{t^{-\frac{n}{2}}} 
{(4\pi)^{\frac{n}{2}}} \tilde{G}_{2,1,n}(\tilde{z}),\  \ \tilde{z} = \frac{|\mathrm{x}|}{2t^{\frac{1}{2}}}
\end{equation}
and then the formula
\begin{equation}
\label{G_G_3}
\tilde{G}_{2,1,n}\left(\frac{z}{\tau}\right) = 2 (4\pi)^{\frac{n}{2}} t^{\frac{\beta n}{\alpha}}\tau^n G_{2,1,n}(\mathrm{x},t^{\frac{2 \beta}{\alpha}}\tau^2),\ z = \frac{|\mathrm{x}|}{2t^{\frac{\beta}{\alpha}}}
\end{equation}
that connects the function $\tilde{G}_{2,1,n}$ from \eqref{sub_3} and the fundamental solution 
$G_{2,1,n}$. Now we put \eqref{G_G_3} into \eqref{sub_3} and get the integral representation
\begin{equation}
\label{sub_37} 
G_{\alpha,\beta,n}(\mathrm{x},t) = \frac{2}{\alpha}\int_0^\infty \tau^{n-1}\, \Psi_{\alpha,\beta,n}(\tau) G_{2,1,n}(\mathrm{x},t^{\frac{2 \beta}{\alpha}}\tau^2)\, d\tau
\end{equation}
that can be rewritten in the form (see \eqref{sub_1_1_3}) 
$$
G_{\alpha,\beta,n}(\mathrm{x},t) = \int_0^\infty t^{-\frac{2\beta}{\alpha}}\Phi_{\alpha,\beta}(st^{-\frac{2\beta}{\alpha}})\, G_{2,1,n}(\mathrm{x},s)\, ds
$$ 
with
\begin{equation}
\label{P-P} 
\Phi_{\alpha,\beta}(\tau) = \frac{1}{\alpha} \tau^{\frac{n}{2}-1} \Psi_{\alpha,\beta,n}\left(\tau^{\frac{1}{2}}\right)
\end{equation}
after the variables substitution $s=t^{\frac{2 \beta}{\alpha}}\tau^2$. 

To determine the kernel function $\Phi_{\alpha,\beta}$ defined by \eqref{P-P}, we first calculate its Mellin integral transform based on the known Mellin integral transform \eqref{Phi_3} of the function $\Psi_{\alpha,\beta,n}$ and the operational relations  \eqref{(1.52)}-\eqref{(1.53)}:
\begin{equation}
\label{Mel-Phi} 
\Phi_{\alpha,\beta}^*(s) = \frac{2}{\alpha} \frac{\Gamma\left( \frac{2}{\alpha}- \frac{2}{\alpha}s\right) \Gamma\left( 1-\frac{2}{\alpha}+ \frac{2}{\alpha}s\right)}
 {\Gamma\left( 1-\frac{2\beta}{\alpha} + \frac{2\beta}{\alpha}s\right) \Gamma\left( 1-s\right)}.
\end{equation}
Thus the function $\Phi_{\alpha,\beta}$ does not depend on the dimension $n$ and can be represented as the Mellin-Barnes integral (inverse Mellin transform of $\Phi_{\alpha,\beta}^*(s)$) as follows:
\begin{equation}
 \label{MB-Phi}
 \Phi_{\alpha,\beta}(\tau) = \frac{2}{\alpha} \frac{1}{2\pi i}\int_{\gamma-i\infty}^{\gamma+i\infty} \frac{\Gamma\left( \frac{2}{\alpha}- \frac{2}{\alpha}s\right) \Gamma\left( 1-\frac{2}{\alpha}+ \frac{2}{\alpha}s\right)}
 {\Gamma\left( 1-\frac{2\beta}{\alpha} + \frac{2\beta}{\alpha}s\right) \Gamma\left( 1-s\right)}\, \tau^{-s}\, ds.
\end{equation}
Now we see that it is the same Mellin-Barnes integral that we dealt with in Example 2 of the previous section (see the formula \eqref{Ex3_1}) and thus the kernel function $\Phi_{\alpha,\beta}(\tau)$ is a probability density function. 

To complete the proof of the theorem, let us now deduce the representation \eqref{Phi_all} of the probability density function $\Phi_{\alpha,\beta}(\tau)$ defined by the Mellin-Barnes integral  \eqref{MB-Phi}.

The general theory of the Mellin-Barnes integrals (see e.g. \cite{Mar83}) says that the integral \eqref{MB-Phi} has three different series representations depending on the relation between the parameters $\alpha$ and $\beta$ and one has to distinguish between three cases: 
$$
\mbox{(i)} \ \beta < \frac{\alpha}{2},\ \  \mbox{(ii)} \ \beta > \frac{\alpha}{2}, \ \ \mbox{and (iii)} \ \beta = \frac{\alpha}{2}.
$$
The reason for this situation is that the integration contour in the Mellin-Barnes integral \eqref{MB-Phi}  can be transformed either to the
loop $L_{-\infty}$ starting and ending at $-\infty$ and encircling
all poles of the function $\Gamma\left( 1-\frac{2}{\alpha}+ \frac{2}{\alpha}s\right)$ (case (i)) or to the
loop $L_{+\infty}$ starting and ending at $+\infty$ and encircling
all poles of the function $\Gamma\left( \frac{2}{\alpha}- \frac{2}{\alpha}s\right)$ (case (ii)) or  to the
loop $L_{-\infty}$  for $0 <\tau<1$ and to the
loop $L_{+\infty}$  for $\tau>1$ (case (iii)). Then the integrals with the integration contours $L_{-\infty}$  or $L_{+\infty}$ can be represented as some series of the hypergeometric type by using the Jordan lemma and the Cauchy residue theorem (see examples already presented above). Now let us shorty discuss the cases (i)-(iii).

\vspace{0.2cm}

\noindent
{\bf Case (i)}: $\beta < \frac{\alpha}{2}$. 

The poles of $\Gamma\left( 1-\frac{2}{\alpha}+ \frac{2}{\alpha}s\right)$ are at the points $s_k = 1-\frac{\alpha}{2} - \frac{\alpha}{2}k,\ k\in \mathbb{N}_0$. The series representation of the Mellin-Barnes integral \eqref{MB-Phi} thus takes the form:
\begin{equation}
\label{MB-Phi-1}
 \Phi_{\alpha,\beta}(\tau) = \frac{2}{\alpha} \sum_{k=0}^\infty
\frac{\alpha}{2} \frac{(-1)^k}{k!} \frac{\Gamma\left( k+1 \right)}
 {\Gamma\left( 1-\beta-\beta k\right) \Gamma\left( \frac{\alpha}{2}+ 
\frac{\alpha}{2} k \right)}\, \tau^{\frac{\alpha}{2}-1+\frac{\alpha}{2}k}.
\end{equation}
Because $\Gamma(k+1) = k!$, the series \eqref{MB-Phi-1} can be expressed in terms of the four parameters Wright function \eqref{wright}
\begin{equation}
\label{MB-Phi-11}
 \Phi_{\alpha,\beta}(\tau) = \tau^{\frac{\alpha}{2}-1}\, W_{(1-\beta,-\beta),(\frac{\alpha}{2},\frac{\alpha}{2})}\left(-\tau^{\frac{\alpha}{2}}\right)
\end{equation}
and we obtained the first part of the formula \eqref{Phi_all}.

\vspace{0.2cm}

\noindent
{\bf Case (ii)}: $\beta > \frac{\alpha}{2}$. 
 
Now we have to take into consideration the poles of $\Gamma\left( \frac{2}{\alpha}- \frac{2}{\alpha}s\right)$ that are located at the points $s_k = 1 + \frac{\alpha}{2}k,\ k\in \mathbb{N}_0$. The series representation of the Mellin-Barnes integral \eqref{MB-Phi} is as follows:
\begin{equation}
\label{MB-Phi-2}
 \Phi_{\alpha,\beta}(\tau) = \frac{2}{\alpha} \sum_{k=0}^\infty
\frac{\alpha}{2} \frac{(-1)^k}{k!} \frac{\Gamma\left( k+1 \right)}
 {\Gamma\left( 1+\beta k \right) \Gamma\left( 
-\frac{\alpha}{2} k \right)}\, \tau^{-1 - \frac{\alpha}{2}k}.
\end{equation}
Because $\Gamma(k+1) = k!$, the series \eqref{MB-Phi-2} can be expressed in terms of the four parameters Wright function \eqref{wright}
\begin{equation}
\label{MB-Phi-21}
 \Phi_{\alpha,\beta}(\tau) = \tau^{-1}\, W_{(1,\beta),(0,-\frac{\alpha}{2})}\left(-\tau^{-\frac{\alpha}{2}}\right)
\end{equation}
and we obtained the second part of the formula \eqref{Phi_all}.

\vspace{0.2cm}

\noindent
{\bf Case (iii)}: $\beta = \frac{\alpha}{2}$. 

In this case, we repeat the calculations made for the case (i) when $0<\tau <1$ and for the case (ii) when $\tau >1$ and apply the reflection formula for the Gamma-function 
$$
\frac{1}{\Gamma(z)\Gamma(1-z)} = \frac{\sin(\pi z)}{\pi},\ z\in \mathbb{C}
$$
to get the last part of the formula \eqref{Phi_all}.

In the rest of the paper, we collect some remarks regarding the subordination formula presented in Theorem \ref{the1} that in our opinion are worth mentioning. 

\begin{Rem}
The kernel function $\Phi_{\alpha,\beta}$ given by \eqref{Phi_all} can be represented as the following Mellin-Barnes integral (see the formula \eqref{MB-Phi} in the proof of Theorem \ref{the1}) for all values of $\alpha$ and $\beta$ under the conditions stated in Theorem \ref{the1}: 
$$
 \Phi_{\alpha,\beta}(\tau) = \frac{2}{\alpha} \frac{1}{2\pi i}\int_{\gamma-i\infty}^{\gamma+i\infty} \frac{\Gamma\left( \frac{2}{\alpha}- \frac{2}{\alpha}s\right) \Gamma\left( 1-\frac{2}{\alpha}+ \frac{2}{\alpha}s\right)}
 {\Gamma\left( 1-\frac{2\beta}{\alpha} + \frac{2\beta}{\alpha}s\right) \Gamma\left( 1-s\right)}\, \tau^{-s}\, ds.
$$
\end{Rem}

\begin{Rem}
The second line of the formula \eqref{Phi_all} can be rewritten in the form
$$
\tau^{-1}\, W_{(1,\beta),(0,-\frac{\alpha}{2})}\left(-\tau^{-\frac{\alpha}{2}}\right) = 
-\tau^{-1-\frac{\alpha}{2}}\, W_{(1-\beta,\beta),(-\frac{\alpha}{2},-\frac{\alpha}{2})}\left(-\tau^{-\frac{\alpha}{2}}\right)  
$$
because the first term of the series in \eqref{MB-Phi-2} or \eqref{MB-Phi-21} is equal to zero due to the fact that the Gamma-function has a pole at the point zero. Thus we can move the summation index $k\to k+1$ in \eqref{MB-Phi-2} and get the representation above. In particular, it is now clear that the kernel function $\Phi_{\alpha,\beta}$ is integrable at $+\infty$. 
\end{Rem}

\begin{Rem}
In the case $\frac{\beta}{\alpha} = \frac{1}{2},\ 0<\beta <1$, the kernel function $\Phi_{\alpha,\beta}(\tau)$ defined by the 3rd line of \eqref{Phi_all} has an integrable singularity at the point $\tau = 1$.
\end{Rem}

\begin{Rem}
The relation \eqref{Ex3_27} between the Mellin integral transform of the kernel function $\Phi_{\alpha,\beta}(\tau)$ and the Mittag-Leffler function \eqref{Ex3_3} can be rewritten in terms of the Laplace integral transform (see the formulas \eqref{cm-1}, \eqref{cm-2}) and thus $\Phi_{\alpha,\beta}(\tau)$ can be also interpreted as the inverse Laplace transform of the Mittag-Leffler function $E_\beta (-\lambda^{\frac{\alpha}{2}})$:
\begin{equation}
\label{Phi_lap} 
E_\beta (-\lambda^{\frac{\alpha}{2}}) = \int_0^\infty \Phi_{\alpha,\beta}(\tau)\, e^{-\lambda \tau}\, d\tau. 
\end{equation}
\end{Rem}

\begin{Rem}
For the time-fractional diffusion equation ($\alpha=2$, $0<\beta \le 1$ in the equation \eqref{eq1}) the subordination formula \eqref{sub_1_1_3} with  the kernel function $\Phi_{\alpha,\beta}(\tau)$ given by the 1st line of \eqref{Phi_all} is valid. In this case, the four parameters Wright function is reduced to the conventional Wright function and we arrive at the known formula (see \eqref{sub_1_1_1})
\begin{equation}
\label{Phi_time} 
G_{2,\beta,n}(\mathrm{x},t) = \int_0^\infty t^{-\beta} W_{1-\beta,-\beta}(-st^{-\beta}) G_{2,1,n}(\mathrm{x},s)\, ds,\  0<\beta < 1.
\end{equation}
\end{Rem}

\begin{Rem}
For the space-fractional diffusion equation ($\beta=1$, $0<\alpha \le 2$ in the equation \eqref{eq1}) the subordination formula \eqref{sub_1_1_3} has to be applied with  the kernel function $\Phi_{\alpha,\beta}(\tau)$ given by the 2nd line of \eqref{Phi_all}. The four parameters Wright function from \eqref{Phi_all} is  reduced to the conventional Wright function and we arrive at the subordination formula 
\begin{equation}
\label{Phi_sp} 
G_{\alpha,1,n}(\mathrm{x},t) = \int_0^\infty s^{-1} W_{0,-\frac{\alpha}{2}}(-s^{-\frac{\alpha}{2}}t) G_{2,1,n}(\mathrm{x},s)\, ds,\ 0<\alpha < 2. 
\end{equation}
\end{Rem}


\begin{thebibliography}{99}

\bibitem{Baz00}
E. Bazhlekova, \textit{Subordination principle for fractional evolution equations}, Fract. Calc. Appl. Anal. \textbf{3} (2000), 213--230.

\bibitem{Baz01}
E. Bajlekova, \textit{Fractional Evolution Equations in Banach Spaces}, Ph.D. thesis, Eindhoven, The Netherlands, 2001.

\bibitem{Baz15}
E. Bazhlekova, \textit{Completely monotone functions and some classes of fractional evolution equations}, Integral Transforms and Special Functions \textbf{26} (2015), no. 9, 737--752.

\bibitem{Baz17}
E. Bazhlekova, I.B. Bazhlekov,  
\textit{Subordination approach to multi-term time-fractional 
diffusion-wave equations}, Journal of Computational and Applied Mathematics, in press, doi: 10.1016/j.cam.2017.11.003.

\bibitem{BL1} L. Boyadjiev, Yu. Luchko,
\textit{Mellin integral transform approach to analyze the multidimensional diffusion-wave equations}, Chaos, Solitons \& Fractals {\bf 102}  (2017), 127--134. 

\bibitem{BL2} L. Boyadjiev, Yu. Luchko,
\textit{Multi-dimensional $\alpha$-fractional diffusion-wave equation and some properties of its fundamental solution}, Computers \& Mathematics with Applications,  {\bf 73} (2017), 2561--2572.

\bibitem{EK2004}
S.D. Eidelman,  A.N. Kochubei, \textit{Cauchy problem for fractional
diffusion equations}, J. Differential Equations,  {\bf 199} (2004), 211--255.

\bibitem{E1953} 
A. Erd\'elyi, \textit{ Higher Transcendental Functions}, vol.2, "McGraw-Hill", New York, 1953.

\bibitem{E1955} 
A. Erd\'elyi, \textit{Higher Transcendental Functions}, vol.3, "McGraw-Hill", New York, 1955.

\bibitem{feller}
W. Feller, 
\textit{ An Introduction to Probability Theory and its Applications},
vol. 2, 2nd edition, "J. Wiley \& Sons Inc.", New York, 1971.

\bibitem{FV2016}
M. Ferreira, N. Vieira, \textit{Fundamental solutions of the time fractional diffusion-wave and parabolic Dirac operators}, 
Journal of Mathematical Analysis and Applications, {\bf 447} (2016),   
329--353.

\bibitem{Fox61}
C. Fox, \textit{The $G$-  and $H$-functions  as  symmetrical  Fourier
kernels}, Trans. Amer. Math. Soc. {\bf 98} (1961), 395--429.

\bibitem{Han2}
A. Hanyga, 
\textit{Multi-dimensional solutions of space-time-fractional diffusion equations}, 
Proc.
R. Soc. Lond. A,   {\bf 458} (2002),  429--450.

\bibitem{Han3}
A. Hanyga, 
\textit{Multidimensional solutions of time-fractional diffusion-wave equations},  Proc.
R. Soc. Lond. A   {\bf 458}  (2002),   933--957.

\bibitem{Kiryakova}
V. Kiryakova, \textit{ Generalized Fractional Calculus and Applications},
"Longman", Harlow, 1994.

\bibitem{K1990}
A.N. Kochubei, \textit{Fractional-order diffusion},  Differential Equations {\bf 26} (1990),  485--492.

\bibitem{K2014}
A.N. Kochubei, \textit{Cauchy problem for fractional diffusion-wave equations with variable coefficients}, Applicable Analysis  {\bf 93} (2014), 2211--2242.

\bibitem{L1999} 
Yu. Luchko,
\textit{Operational method in fractional calculus},  Fract. Calc. Appl. Anal., {\bf 2} (1999), 463--489.

\bibitem{L2010} 
Yu. Luchko,
\textit{Some uniqueness and existence results for the initial-boundary-value problems for the generalized time-fractional diffusion equation}, Computers and Mathematics with Applications, {\bf 59} (2010), 1766--1772.

\bibitem{L2013} 
Yu. Luchko,
\textit{Fractional wave equation and damped waves},   J. Math. Phys., {\bf 54} (2013), 031505.

\bibitem{L2014} 
Yu. Luchko,
\textit{Multi-dimensional fractional wave equation and some properties of its fundamental solution}, Communications in Applied and Industrial Mathematics, {\bf 6}  (2014), e-485. 

\bibitem{L2015} 
Yu. Luchko,
\textit{Wave-diffusion dualism of the neutral-fractional processes}, Journal of Computational Physics, {\bf 293} (2015), 40--52.

\bibitem{L2016} 
Yu. Luchko,
\textit{A new fractional calculus model for the two-dimensional anomalous diffusion and its analysis}, Math. Model. Nat. Phenom., {\bf 11} (2016), 1--17.

\bibitem{L2017} 
Yu. Luchko,
\textit{On some new properties of the fundamental solution to the multi-dimensional space- and time-fractional diffusion-wave equation}, Mathematics, {\bf 5} (2017), no. 4, 1--16. 

\bibitem{Luc}
Yu. Luchko, R. Gorenflo, \textit{Scale-invariant solutions of a partial
differential equation of fractional order}, Fract. Calc. Appl. Anal.  {\bf 1} (1998), 63--78.

\bibitem{LK2013} 
Yu. Luchko, V. Kiryakova, 
\textit{The Mellin integral transform in fractional calculus},  Fract. Calc. Appl. Anal., {\bf 16} (2013), 405--430.

\bibitem{Mai-P} 
F. Mainardi, G. Pagnini,  
\textit{  Salvatore Pincherle: The pioneer of the Mellin-Barnes integrals}, 
J. Computational and Applied  Mathematics, {\bf 153} (2003), 331--342.

\bibitem{MLP2001} 
F. Mainardi, Yu. Luchko, G. Pagnini, \textit{ The fundamental solution of the space-time fractional diffusion equation}, Fract. Calc. Appl. Anal.,  {\bf 4} (2001),  153--192.

\bibitem{Mar83}
O.I. Marichev, \textit{Handbook   of  integral  transforms  of
higher transcendental  functions,  theory  and  algorithmic  tables}, 
"Ellis Horwood", Chichester, 1983. 

\bibitem{MatSax78}
A.M. Mathai, R.K. Saxena, \textit{ The $H$-functions with
 Applications in Statistics and Other Disciplines}, "John Wiley", New
 York, 1978.
 
 \bibitem{Miller-Samko}
K.S. Miller, S.G. Samko, \textit{ Completely monotonic functions}, 
Integral Transforms and Special  Functions, {\bf 12} (2001), 389--402.

\bibitem{Pru93}
 J. Pr\"{u}ss, \textit{Evolutionary Integral Equations and Applications}, 
"Birkh\"{a}user", Basel, 1993.
 
\bibitem{SZ1997} 
A. Saichev, G. Zaslavsky,  \textit{Fractional kinetic equations: Solutions and applications}, Chaos {\bf 7} (1997), 753--764.

\bibitem{SKM1993} 
S.G. Samko, A.A. Kilbas, O.I. Marichev, O.I.  \textit{Fractional Integrals and Derivatives: Theory and Applications}, "Gordon and Breach", New York, 1993.

\bibitem{[SSV]}
R.L. Schilling, R. Song,  Z. Vondra\v cek, \textit{ Bernstein
Functions. Theory and Applications}, "De Gruyter", Berlin, 2010.

\bibitem{SW1989} 
W.R. Schneider, W. Wyss,  \textit{Fractional diffusion and wave equations},   J. Math. Phys.,  {\bf 30} (1989), 134--144.

\bibitem{Wr2}
E.M. Wright, \textit{The asymptotic expansion of the generalized
hypergeometric function}, Journal London Math. Soc., {\bf 10} (1935),  
287--293.

\bibitem{YakLuc94}
S. Yakubovich, Yu. Luchko, \textit{The Hypergeometric Approach to
Integral Transforms and Convolutions}, "Kluwer Acad. Publ.", Dordrecht, 1994.

\end{thebibliography}
\end{document}